\documentclass[]{article}
\usepackage{amsmath,bm}
\numberwithin{equation}{subsection}
 
\usepackage{algorithm}
\usepackage{tikz,pgfplots}
\usepackage{algpseudocode}
\usepackage{amssymb}
\usepackage{graphicx}
\usepackage{color}
\usepackage{cite}
\usepackage{amsthm}
\usepackage{url, verbatim}

\usepackage{chngcntr}
\counterwithout{equation}{subsection}
\counterwithin{equation}{section}

\theoremstyle{remark}

\usepackage[left=1.30in,right=1.30in,top=1.5in]{geometry}

\usepackage{bbm,mathabx}

\usepackage{array,multirow,booktabs} 

\usepackage{amsfonts}
\usepackage{algorithmicx}
\usepackage{algpseudocode}
\usepackage{epstopdf}
\usepackage{mathrsfs,mathtools}
\usepackage{xcolor}
\usepackage{enumitem}
\usepackage{chngcntr}

\makeatletter
\def\BState{\State\hskip-\ALG@thistlm}
\makeatother

\newcommand{\bs}[1]{\boldsymbol{#1}}
\newcommand{\bmu}{\bs{\mu}}

\newcommand{\R}{\mathbbm{R}}

\newcommand{\calC}{\mathcal{C}}

\newcommand{\calO}{\mathcal{O}}
\newcommand{\calD}{\mathcal{D}}

\newcommand{\calL}{\mathcal{L}}
\newcommand{\calN}{\mathcal{N}}

\newcommand{\Nmax}{N_{\mathrm max}}

\DeclareMathOperator*{\argmax}{argmax}

\newcommand{\naive}{na\"{\i}ve}

\begin{document}

\title{Offline-Enhanced Reduced Basis Method through adaptive construction of the Surrogate Parameter Domain 
}
\author{Jiahua Jiang\thanks{
Department of Mathematics, University of Massachusetts Dartmouth, 285 Old Westport Road, North Dartmouth, MA 02747, USA. Emails: {\tt{\{jjiang, yanlai.chen\}@umassd.edu}}. J.~Jiang {was} partially supported by National Science Foundation grant DMS-1216928.
}
\and
Yanlai Chen\footnotemark[1]
\and
Akil Narayan\thanks{
Scientific Computing and Imaging (SCI) Institute and Department of Mathematics, University of Utah, 72 S Campus Drive, Salt Lake City, UT 84112, USA. Email: {\tt{akil@sci.utah.edu}.} A.~Narayan is partially supported by NSF DMS-1552238, AFOSR FA9550-15-1-0467, and DARPA N660011524053
}
}

\date{}
\maketitle

\begin{abstract}
The Reduced Basis Method (RBM) is a popular certified model reduction approach for solving parametrized partial differential equations. One critical stage of the \textit{offline} portion of the algorithm is a greedy algorithm, requiring maximization of an error estimate over parameter space. In practice this maximization is usually performed by replacing the parameter domain continuum with a discrete "training" set. When the dimension of parameter space is large, it is necessary to significantly increase the size of this training set in order to effectively search parameter space. Large training sets diminish the attractiveness of RBM algorithms since this proportionally increases the cost of the offline {phase}.

In this work we propose novel strategies for offline RBM algorithms that mitigate the computational difficulty of maximizing error estimates over a training set. The main idea is to identify a subset of the training set, a "surrogate parameter domain" (SPD), on which to perform greedy algorithms. The SPD's we construct are much smaller in size than the full training set, yet our examples suggest that they are accurate enough to represent the solution manifold of interest at the current offline RBM iteration. We propose two algorithms to construct the SPD: Our first algorithm, the Successive Maximization Method (SMM) method, is inspired by inverse transform sampling for non-standard univariate probability distributions. The second constructs an SPD by identifying pivots in the Cholesky Decomposition of an approximate error correlation matrix. We demonstrate the algorithm through numerical experiments, showing that the algorithm is capable of accelerating offline RBM procedures without degrading accuracy, assuming that the solution manifold has low Kolmogorov width.
\end{abstract}

\section{Introduction}
Increasing computer power and computing availability makes the simulation of large-scale or high-dimensional numerical problems modeling complex phenomena more accessible, but any reduction in computational effort is still required for sufficiently onerous many-query computations and repeated output evaluations for different values of some inputs of interest. Much recent research has concentrated on schemes to accelerate computational tools for such many-query or parameterized problems; these schemes include Proper Orthogonal Decomposition (POD) \cite{GPJ1993, KJ2002}, balanced truncation method \cite{SA2004, B1981}, Krylov Subspace methods \cite{PR1995, Y1996} and the Reduced Basis Method (RBM) \cite{RozzaHuynhPatera2008, MR2821591}.  The basic idea behind each of these model order reduction techniques is to iteratively project an associated large algebraic system to a small system that can effectively capture most of the information carried by the original model. In the context of parameterized partial differential equations (PDE), the fundamental reason that such model reduction approaches are accurate is that, for many PDE's of interest, the solution manifold induced by the parametric variation has small Kolmogorov width \cite{MR774404}. 

Among the strategies listed above, one of most appealing methods is the Reduced Basis Method. RBM seeks to parametrize the random inputs and select the most representative points in the parameter space by means of a greedy algorithm that leverages an \textit{a posteriori} error estimate. RBM algorithm are split into \textit{offline} and \textit{online} stages. During the offline stage, the parameter dependence is examined and the greedy algorithm is used to judiciously select a small number of parameter values on which the full, expensive PDE solver is solved. The solutions on this small parameter set are called \textit{snapshots}. During the online stage, {an} approximate solution for any new parameter value is efficiently computed via a linear combination of the offline-computed snapshots. This linear combination can usually be accomplished with orders of magnitude less effort than a full PDE solve. Thus, RBM invests significant computational effort in an offline stage so that the online stage is efficient \cite{MR804937, MR616719, MR917452, MR2824231, MR3129759}.

The Reduced Basis Method was initially introduced for structure analysis \cite{AJ1980} and recently has undergone vast development in theory \cite{RozzaHuynhPatera2008, YAG2002, AG2007, MA2005, MYNA2007} and applied to many engineering problems \cite{MR2326199, MR2639602, MR2824231, MR2914321, MR3127301, MR3475656}. If the input parameters to the PDE are random, then RBM is philosophically similar to stochastic collocation \cite{MR2646806}, but uses an adaptive sampling criterion dictated by the PDE's \textit{a posteriori} error estimate. It has been shown that RBM can help {delay} the curse of dimensionality when solving parameterized problems in uncertainty quantification whenever the solution manifold lies in a low dimensional space \cite{jiang_goal-oriented_2016}.

RBM is motivated by the observation that the parameteric variation in many parameterized PDEs can be well-approximated by a finite-dimensional projection. Early RBM research concentrated on problems with a low-dimensional parameter due to the lack of effective tools to sample the \textit{a posteriori} error estimate over high-dimensional spaces \cite{RozzaHuynhPatera2008}. An effective procedure for greedily selecting parameter values for use in the RBM procedure must simultaneously leverage the structure of the error estimate along with efficient methods for searching over high-dimensional spaces.
Some recent effort in the RBM framework has been devoted {to} inexpensive computation of the \textit{a posteriori} error estimate along with effective sampling strategies \cite{MR2123791}.

The portion of the RBM algorithm most relevant in the context of this article is in the offline stage: Find a parameter value that maximizes a given (computable) error estimate. This maximization is usually accomplished in the computational setting by replacing the parameter domain continuum by a large, discrete set called the \textit{training set}. Even this \naive{} procedure requires us to compute the value of the error estimate at every point in the training set. Thus, the work required scales proportionally to the training set size. When the parameter is high-dimensional, the size of the training set must be very large if one seeks to search over all regions of parameter space. This onerous cost of the offline stage debilitates RBM in this scenario. Thus, assuming the training set must be large, a more sophisticated scheme for maximizing the error estimate must be employed. To the best of our knowledge, computational stratagems in the RBM framework to tackle this problem are underdeveloped.

In this paper we propose two novel strategies for mitigating the cost of searching over a training set of large size. The essential idea in both approaches is to perform some computational analysis on the \textit{a posteriori} error estimate in order to construct a \textit{surrogate domain}, a subset of the original training set with a much smaller size, that can effectively predict the general trend of the error estimate. The construction of this surrogate domain must be periodically repeated during the iterative phase of the offline RBM algorithm. Ideally, we want to decrease the computational burden of the offline algorithm without lowering the fidelity of the RBM procedure. The following qualitative characteristics are the guiding desiderata for construction of the surrogate domain:
\begin{enumerate}[noitemsep]
  \item the information used for surrogate domain construction should be inexpensive to obtain
  \item the parameteric variation on surrogate domain should be representative of that in the original training set
  \item the size of the surrogate domain should be significantly smaller than that of the original training set
\end{enumerate}
Our proposed offline-enhanced RBM strategies are as follows: 
\begin{enumerate}[noitemsep]
  \item Successive Maxmimization Method for Reduced Basis (SMM-RBM) --- We construct an empirical cumulative distribution function of the \textit{a posteriori} error estimate on the training set, and deterministically subsample a surrogate domain according to this distribution.
  \item Cholesky Decomposition Reduced Basis Method (CD-RBM) --- An approximate correlation matrix (Gramian) of errors over the training set is computed, and the pivots in a pivoted Cholesky decomposition \cite{MR2899254} identify the surrogate domain. 
\end{enumerate}
We note that both of our strategies are empirical in nature. In particular, it is relatively easy to manufacture error estimate data so that our construction of a surrogate domain does not accurately capture the parameter variation over the full training set. One of the main observations we make in our numerical results section is that such an adversarial situation does not occur for the parameterized PDEs that we investigate. Our procedure also features some robustness: a poorly constructed surrogate domain does not adversely affect either the efficiency or the accuracy of the RBM simulation.

The remainder of this paper is organized as follows. A parametrized PDE with random input data is set up with appropriate assumptions on the PDE operator in section \ref{sec:background}. The general framework and properties of the Reduced Basis Method are likewise introduced in Section \ref{sec:background}. Section \ref{sec:offline-rbm} is devoted to the development of our novel offline-enhanced RBM methods, consisting of SMM-RBM and CD-RBM. A rough complexity analysis of these methods is given in section \ref{sec:oerbm-complexity}. Our numerical examples are shown in section \ref{sec:results}. 

\section{Background}\label{sec:background}
In this section, we introduce the Reduced Basis Method (RBM) in its classical form; much of this is standard in the RBM literature. The reader familiar with RBM methodology may skip this section, using Table \ref{tab:notation} as a reference for our notation.

\begin{table}
  \begin{center}
  \resizebox{\textwidth}{!}{
    \renewcommand{\tabcolsep}{0.4cm}
    \renewcommand{\arraystretch}{1.3}
    {\scriptsize
    \begin{tabular}{@{}lp{0.8\textwidth}@{}}
      \toprule
      $\bmu$ & Parameter in $\calD \subseteq \R^p$ \\
      $u(\bmu)$ & Function-valued solution of a parameterized PDE \\
      $\mathcal{N}$ & Degrees of freedom (DoF) in PDE ``truth" solver \\
      $u^{\mathcal{N}}(\bmu)$ & Truth solution (finite-dimensional)\\
      $N$ & Number of reduced basis snapshots, $N \ll \mathcal{N}$\\
      $\bmu^j$ & ``Snapshot" parameter values, $j=1, \ldots, N$\\
      $X^{\mathcal{N}}_{{N}}$ & Span of $u^{\mathcal{N}}\left(\bmu^k\right)$ for $k=1, \ldots, N$\\
      $u_N^{\mathcal{N}}(\bmu)$ & Reduced basis solution, $u_N^{\mathcal{N}} \in X^{\mathcal{N}}_{{N}}$\\
      $e_N(\bmu)$ & Reduced basis solution error, equals $u^{\mathcal{N}}(\bmu) - u_N^{\mathcal{N}}(\bmu)$ \\
      $\Xi_{\rm{train}}$ & Parameter training set, a finite subset of $\mathcal{D}$ \\
      $\Delta_{{N}} \left(\bmu\right)$ & Error estimate (upper bound) for $\left\|e_N\left(\bmu\right)\right\|$ \\
      $\epsilon_{\mathrm{tol}}$ & Error estimate stopping tolerance in greedy sweep \\
    \bottomrule
    \end{tabular}
  }
    \renewcommand{\arraystretch}{1}
    \renewcommand{\tabcolsep}{12pt}
  }
  \end{center}
  \caption{Notation used throughout this article.}\label{tab:notation}
\end{table}

\subsection{Problem setting}
Let $\mathcal{D} \subset \mathbb{R}^{p}$ be the range of variation of a $p$-dimensional parameter and $\Omega \subset \mathbb{R}^{d}$ (for $d = 2 ~ \text{or} ~ 3$) a bounded spatial domain. We consider the following parametrized problem: Given $\bm{\mu} \in \mathcal{D}$, the goal is to evaluate the output of interest
\begin{equation}
\label{eq:interest}
s(\bm{\mu}) = \ell(u(\bm{\mu})){,}
\end{equation}
where the function $u(\bm{\mu}) \in X$ satisfies
\begin{equation}\label{eq:sat}
a(u(\bm{\mu}),v,\bm{\mu}) = f(v; \bm{\mu}), \quad v \in X,
\end{equation}
which is a parametric partial differential equation (pPDE) written in a weak form. Here $X = X(\Omega)$ is a Hilbert space satisfying $H^{1}_{0}(\Omega) \subset X(\Omega) \subset H^{1}(\Omega)$. We denote by $(\cdot , \cdot)_{X}$ the inner product associated with the space $X$, whose induced norm $|| \cdot ||_{X} = \sqrt{(\cdot , \cdot)_{X}}$ is equivalent to the usual $H^{1}(\Omega)$ norm. 
We assume that $a(\cdot ,\cdot; \bm{\mu}): X \times X \rightarrow \mathbb{R}$ is continuous and coercive over $X$ for all $\bm{\mu}$ in $\mathcal{D}$, that is,
\begin{subequations}
\begin{equation}\label{eq:con}
\gamma(\bm{\mu}) \coloneqq \underset{w \in X}{\sup} ~ \underset{v \in X}{\sup}~ \frac{a(w, v; \bm{\mu})}{||w||_{X}||v||_{X}} < \infty, \quad \forall \bm{\mu} \in \mathcal{D}, 
\end{equation} 
\begin{equation}\label{eq:coe}
\alpha(\bm{\mu}) \coloneqq \underset{w \in X} \inf \frac{a(w, w; \bm{\mu})}{||w||^2_{X}} \geqslant \alpha_{0} > 0, \forall \bm{\mu} \in \mathcal{D}. 
\end{equation}
\end{subequations}
$f(\cdot)$ and $\ell(\cdot)$ are linear continuous functionals over $X$, and for simplicity we assume that $\ell$ is independent of $\bm{\mu}$. 

We assume that $a(\cdot, \cdot; \bm{\mu})$ is ``affine'' with respect to functions of the parameter $\bm{\mu}$: there exist $\bm{\mu}$-dependent coefficient functions $\Theta_{a}^{q}: \calD \rightarrow \mathbb{R}$ for $q = 1, \ldots Q_a$, and corresponding continuous $\bmu$-independent bilinear forms $a^{q}(\cdot, \cdot): X \times X \rightarrow \R$ such that 
\begin{equation}\label{eq:assum_a}
a(w, v; \bm{\mu}) = \sum_{q = 1}^{Q_{a}} \Theta_{a}^{q}(\bm{\mu})a^{q}(w,v).
\end{equation}
This assumption of affine parameter dependence is common in the reduced basis literature \cite{RozzaHuynhPatera2008}, and remedies are available \cite{BarraultMaday2004} when it is not satisfied.

Finally, we assume that there is a finite-dimensional discretization for the model problem \eqref{eq:sat}: The solution space $X$ is discretized by an $\calN$-dimensional subspace $X^{\mathcal{N}}$ (i.e., $dim(X^{\mathcal{N}}) = \mathcal{N}$) and \eqref{eq:interest} and \eqref{eq:sat} are discretized as
\begin{equation}\label{eq:update_problem}
\begin{cases}
\text{For } \bm{\mu} \in \mathcal{D},~ \text{solve} \\
s^{\mathcal{N}} = \ell(u^{\mathcal{N}}(\bm{\mu})) ~ \text{where} ~ u^{\mathcal{N}}(\bm{\mu}) \in X^{\mathcal{N}} ~ \text{satisfies} \\
a(u^{\mathcal{N}}, v; \bm{\mu}) = f(v; \bm{\mu}) \quad \forall v \in X^{\mathcal{N}}.
\end{cases}
\end{equation}
The relevant quantities such as the coercivity constant \eqref{eq:coe} are defined according to the discretization,
\[
\alpha^{\mathcal{N}}(\bm{\mu})  = \inf_{w \in X^{\mathcal{N}}} \frac{a(w, w; \bm{\mu})}{||w||^{2}_{X}}, ~ \forall \bm{\mu} \in \mathcal{D}.
\]
In the RBM literature, any discretization associated to $\mathcal{N}$ is called a ``truth" discretization. E.g., $u^{\mathcal{N}}$ is called the ``truth solution".

\subsection{RBM framework}
We assume $\calN$ is large enough so that solving \eqref{eq:update_problem} gives highly accurate approximations for $\bmu \in \mathcal{D}$. However, large $\calN$ also means that solving \eqref{eq:update_problem} is expensive, and in many-query contexts (e.g., optimization) a direct approach to solving \eqref{eq:update_problem} is computationally infeasible. The situation is exacerbated when the $\bmu \mapsto u^{\calN}(\bmu)$ response is sought in a real-time fashion. The reduced basis method is a reliable model reduction tool for these scenarios. 

This section presents a brief overview of the standard RBM algorithm. Given a finite training set of parameter samples $\Xi_{\rm train} \subset \mathcal{D}$ as well as a prescribed maximum dimension $N_{\rm max}$ (usually $\ll \mathcal{N}$), we approximate the solution set
\begin{align*}
  \left\{ u(\bmu) \;\; |\;\; \bmu \in \mathcal{D} \right\} \subset X^{\mathcal{N}}
\end{align*}
via an $N$-dimensional subspace of $X^{\mathcal{N}}$, with $N \leq N_{\rm max}$. In RBM, this is accomplished via the $N$-dimensional cardinal Lagrange space
\begin{equation}\label{eq:space_define}
X^{\mathcal{N}}_{N} \coloneqq \text{span} \{ u^{\mathcal{N}}(\bm{\mu}^{n}), 1 \leq n \leq N \}, \quad N = 1, \dots, N_{\rm max}
\end{equation}
in a hierarchical manner by iteratively choosing samples $S_{N} = \{\bm{\mu}^{1}, \dots, \bm{\mu}^{N}\}$ from the training set $\Xi_{\rm train}$ until $N$ is large enough so that a prescribed accuracy tolerance $\epsilon_{\mathrm tol}$ is met. 
The $u^{\mathcal{N}}(\bm{\mu}^{n})$ for $1 \le n \le N_{\rm max}$ are the so-called ``snapshots'', and are obtained by solving \eqref{eq:update_problem} with $\bm{\mu} = \bm{\mu}^{n}$. 

It is obvious that both $S_{N}$ and $X^{\mathcal{N}}_{N}$ are nested; that is, $S_{1} \subset S_{2} \subset \dots \subset S_{N_{\rm max}}$ and $X^{\mathcal{N}}_{1} \subset X^{\mathcal{N}}_{2} \subset \dots \subset X^{\mathcal{N}}_{N_{\rm max}} \subset X^{\mathcal{N}}$. This condition is fundamental in ensuring efficiency of the resulting RB approximation. 
Given $\bm{\mu} \in \mathcal{D}$, we seek a surrogate RB solution $u_{N}^{\mathcal{N}}(\bmu)$ in the reduced basis space $X^{\mathcal{N}}_{N}$ for the truth approximation $u^\calN(\bm{\mu})$ by solving the following {\em reduced} system
\begin{equation}
\label{eq:reduced_system}
\begin{cases}
\text{For } \bm{\mu} \in \mathcal{D},~ \text{evaluate} \\
s_{N}^{\mathcal{N}} = \ell(u_{N}^{\mathcal{N}}(\bm{\mu})) ~ \text{s.t.} ~ u_{N}^{\mathcal{N}}(\bm{\mu}) \in X_{N}^{\mathcal{N}} \subset X^{\mathcal{N}}~ \text{satisfies} \\
a(u_{N}^{\mathcal{N}}, v; \bm{\mu}) = f(v) \quad \forall v \in X_{N}^{\mathcal{N}}.
\end{cases}
\end{equation}
In comparison to the $\mathcal{N}$-dimensional system \eqref{eq:update_problem}, the reduced system \eqref{eq:reduced_system} is $N$-dimensional; when $N \ll \mathcal{N}$, this results in a significant computational savings. The Galerkin procedure in \eqref{eq:reduced_system} selects the best solution in $X_N^{\calN}$ satisfying the pPDE.
\footnote{In implementations, in order to ameliorate ill-conditioning issues that may arise in \eqref{eq:reduced_system} we first apply the Gram-Schmidt process with respect to the $(\cdot , \cdot)_{X}$ inner product each time a new snapshot $u^{\mathcal{N}}(\bm{\mu}^{n})$ is generated to obtain a $(\cdot , \cdot)_{X}$-orthonormal basis $\{\xi_{n}^{\mathcal{N}}\}_{n = 1}^{N_{\rm max}}$. We omit explicitly denoting or showing this orthogonalization procedure.} Then, the RB solution $u_{N}^{\mathcal{N}}(\bm{\mu})$ for any parameter $\bm{\mu} \in \mathcal{D}$ can be expressed as  
\begin{equation}\label{eq:rbsolution}
  u_{N}^{\mathcal{N}}(\bm{\mu}) = \sum_{m = 1}^{N} u_{Nm}^{\mathcal{N}}(\bm{\mu}) u^{\mathcal{N}}\left(\bmu^n\right) 
\end{equation}
Here $\{u_{Nm}^{\mathcal{N}}(\bm{\mu})\}_{m = 1}^{N}$ are the unknown RB coefficients that {can} be obtained by solving \eqref{eq:reduced_system}. Upon replacing the reduced basis solution in \eqref{eq:reduced_system} by \eqref{eq:rbsolution} and taking the $X_N^{\mathcal{N}}$ basis functions $v_n = u^{\mathcal{N}} \left( \bmu^n\right)$, $1 \le n \le N$, as the test functions for Galerkin approximation, we obtain the RB ``stiffness'' equations
\begin{equation}\label{eq:rbstiff}
\sum_{m = 1}^{N} a(v_m , v_n ; \bm{\mu})u_{Nm}^{\mathcal{N}}(\bm{\mu}) = f(v_n;\bmu) \quad 1\le {n} \le N
\end{equation}
Once this system is solved for the coefficients $u_{Nm}(\bmu)$, the RB output $s_{N}^{\mathcal{N}}(\bm{\mu})$ can be subsequently evaluated as
\begin{equation}\label{eq:rb_out}
s_{N}^{\mathcal{N}}(\bm{\mu}) = \ell(u^\calN_N(\bmu)).
\end{equation}

It is not surprising that the accuracy of the RB solution $u_{N}^{\mathcal{N}}$ and of the corresponding computed output of interest $s_{N}^{\mathcal{N}}$ both depend crucially on the construction of the reduced basis approximation space. The procedure we use for efficiently selecting representative parameters $\bmu^1, \ldots, \bmu^N$ and the corresponding snapshots defining the reduced basis space plays an essential role in the reduced basis method. 

\subsection{Selecting snapshots: Enhancing offline RBM computations}
This paper's main contribution is the development of novel procedures for selecting the snapshot set $\left\{\bmu^n\right\}_{n=1}^{\Nmax}$. The main idea of our procedure is very similar to classical RBM methods, the latter of which is the greedy scheme
\begin{align}\label{eq:greedy-selection}
  \bmu^{n+1} &= \argmax_{\bmu \in \Xi_{\rm train}} \Delta_n(\bmu), & \Delta_n(\mu) \geq \left\| u^\calN_n(\bmu) - u^\calN(\bmu)\right\|_{X^{\calN}}.
\end{align}
Here $\Delta_n(\cdot)$ is an efficiently-computable error estimate for the RBM procedure, and $\Xi_{\rm train} \subset \calD$ is a \textit{training set}, i.e., a large but finite discrete set that replaces the continuum $\calD$. Once $\mu^{n+1}$ is selected, standard RBM mechanics can be used to construct $u^{\calN}_{n+1}$ so that the procedure above can be iterated to compute $\bmu^{n+2}$.

We present more details of this procedure in the Appendix, including a mathematical justification of why the above procedure is effective and computable. A classical RBM algorithm computes the maximum over $\Xi_{\rm train}$ above in a brute-force manner; since $\Xi_{\rm train}$ is large and this maximization must be done for every $n =1, \ldots, \Nmax$, this process of selecting snapshots is usually one of the more computationally expensive portions of RBM algorithms. 

This manuscript is chiefly concerned with ameliorating the cost of selecting snapshots; we call this an ``offline-enhanced" {Reduced Basis Method}\footnote{``Offline" is a standard descriptor for this general portion of the full RBM algorithm; see the appendix.}. We present two algorithms that are alternatives to the brute-force approach \eqref{eq:greedy-selection}. Instead of maximizing $\Delta_n$ over the full training set $\Xi_{\rm train}$, we instead maximize over a subset of $\Xi_{\rm train}$ that we call the ``surrogate domain". The efficient computational determination of the surrogate domain, and subsequent empirical studies investigating the accuracy and efficiency of offline-enhanced methods compared to classical {RBM}, are the remaining topics of this manuscript.

This paper does not make novel contributions to any of the other important aspects of RBM algorithms (offline/online decompositions, error estimate computations, etc.), but in the interests of completeness we include in the Appendix a brief overview of the remaining portions of RBM algorithms.

\section{Offline-enhanced RBM: Design and analysis}\label{sec:offline-rbm}

The complexity of the offline stage, where the optimization \eqref{eq:greedy-selection} is performed, depends on $N_{\rm train}$. Although this dependence is only linear, a large $N_{\rm train}$ can easily make the computation onerous; such a situation arises when the parameter domain $\mathcal{D}$ has large dimension $p$. In this case standard constructions for $N_{\rm train}$ yield training sets that grow exponentially with $p$, even when the more parsimonious sparse grid constructions are involved \cite{xiu_high-order_2005}. 
Our goal in this project is to ameliorate the cost of sweeping over a very large training set in \eqref{eq:greedy-selection} without sacrificing the quality of the reduced basis solution. We call this approach an \textit{Offline-enhanced Reduced Basis Method}.

\begin{table}
  \begin{center}
  \resizebox{\textwidth}{!}{
    \renewcommand{\tabcolsep}{0.4cm}
    \renewcommand{\arraystretch}{1.3}
    {\scriptsize
    \begin{tabular}{@{}lp{0.8\textwidth}@{}}
      \toprule
          & {\textbf{Algorithm Template Notation}} \\
    $\ell$ & Number of ``outer" loops in the offline enhancement procedure Algorithm \ref{alg:oerbm} \\
      $E_{\ell}$ & The largest error estimator at the beginning of outer loop $\ell$\\
      $N_{\ell}$ & Number of offline-enhanced snapshots chosen at iteration $\ell$\\
      $\Xi_{\mathrm{sur}}$ & ``Surrogate" parameter domain (SPD), a subset of $\Xi_{\rm{train}}$  \\
      $M_\ell$ & SPD size at outer loop iteration $\ell$ \\
      $K_{\mathrm{damp}}$ & A constant integer controlling the damping ratio $\frac{1}{K_{\rm damp} \times (\ell + 1)} \in (0,1)$ \\
      \midrule
          & \textbf{SMM Notation} \\
      $I^{M}_k$ & $M$ equispaced samples on the interval $(\epsilon_{\mathrm{tol}}, \max_{\bmu \in \Xi_{\rm{train}}} \Delta_k(\bmu)]$ \\
      \midrule
          & \textbf{CD Notation} \\
      $G$ & $\Xi_{\rm{train}} \times \Xi_{\rm{train}}$ Gramian error matrix, with entries $G_{i,j} = \left( e_N\left(\bmu^i\right), e_N\left(\bmu^j\right) \right)_{X^{\calN}}$ \\
      $\widetilde{G}$ & $\Xi_{\rm{train}} \times \Xi_{\rm{train}}$ approximate Gramain error matrix, with entries $\widetilde{G}_{i,j} = \left( \widetilde{e}_N\left(\bmu^i\right), \widetilde{e}_N\left(\bmu^j\right) \right)_{X^{\calN}}$ \\
    \bottomrule
    \end{tabular}
  }
    \renewcommand{\arraystretch}{1}
    \renewcommand{\tabcolsep}{12pt}
  }
  \end{center}
  \caption{Notation used for offline-enhanced RBM.}\label{tab:notation_oeRBM}
\end{table}

The basic idea of our approach is to perform the standard RBM greedy algorithm on a ``surrogate" parameter domain (SPD) constructed as subsets of the original training set $\Xi_{\rm train}$. The SPD is constructed adaptively, and construction is periodically repeated after a small batch of snapshots are selected. We let $\Xi_{\rm Sur}$ denote these constructed SPD's; they are small enough compared to $\Xi_{\rm train}$ to offer considerable acceleration of the greedy sweep \eqref{eq:greedy-selection}, yet large enough to capture the general landscape of the solution manifold. We present in Algorithm \ref{alg:oerbm} a general template for our Offline-Enhanced Reduced Basis Method. This algorithm can be implemented once we describe how $\Xi_{\rm Sur}$ are constructed; these descriptions are the topic of the next sections. In Table \ref{tab:notation_oeRBM} we summarize the notation in Algorithm \ref{alg:oerbm}.

The first contribution of our paper resides in the unique structure of this template. Each global greedy sweep (i.e., over $\Xi_{\rm train}$, and labeled ``One-step greedy" in Algorithm \ref{alg:oerbm}) is followed by multiple targeted sweeps over the (smaller) SPD $\Xi_{\rm Sur}$ (labeled ``Multi-step greedy"). These latter sweeps produce a computational savings ratio of $1 - |\Xi_{\rm Sur}|/|\Xi_{\rm train}|$ because they operate on $\Xi_{\rm Sur}$ instead of on $\Xi_{\rm train}$.

Note that we still require occasional global greedy sweeps, even though they are expensive. These global sweeps are necessary to retain reliability of the greedy algorithm.

\begin{algorithm}[H] 
\begin{algorithmic}[1]
\State Input: training set $\Xi_{\rm train}$, an accuracy tolerance $\varepsilon_{\mathrm{tol}}$.
\State Randomly select the first sample $\bmu^1 \in \Xi_{\rm train}$, and set $n = 1$, $\varepsilon = 2 \varepsilon_{\mathrm{tol}}$, and $\ell = 0$.
\State Obtain truth solution $u^\mathcal{N}(\bmu^1)$, and set $X^\calN_1 = \mbox{span}\left\{u^{\mathcal N}(\bmu^1)\right\}$.
\While {$(\varepsilon > \varepsilon_{\mathrm{tol}})$}
\vspace{0.05in}
\State Set $\ell \leftarrow \ell + 1$.
\State 
\begin{minipage}{0.03\textwidth}
\rotatebox{90}{One-step greedy}
\end{minipage} 
\begin{minipage}{0.03\textwidth}
\rotatebox{90}{scan on $\Xi_{\rm train}$}
\end{minipage} 
\hspace{4mm}
\framebox[0.8 \textwidth]
{\begin{minipage}{0.8 \textwidth}
\For{each $\boldsymbol{\mu} \in \Xi_{\rm train}$}
\State Obtain RBM solution $u^{\mathcal{N}}_{n}(\boldsymbol{\mu}) \in X^\calN_n$ and error estimate ${\Delta_{n}}(\boldsymbol{\mu})$
\EndFor
\vspace{0.05in}
\State $\bmu^{n+1} = \underset{\bm{\mu} \in \Xi_{\rm train}}{\argmax} \Delta_{n}(\bm{\mu})$, $\varepsilon = \Delta_{n}(\bmu^{n+1})$, $E_{\ell} = \varepsilon$.
\State Augment RB space $X^\calN_{n+1} = X^\calN_n \oplus \{u^\calN(\bmu^{n+1})\}$.
\State Set $n \leftarrow n+1$, $N_\ell \gets 0$
\end{minipage}
}
\vspace{0.05in}
\State \hspace{-0.0\textwidth} Construct SPD $\Xi_{\rm Sur}$ based on $\{(u_{n-1}^\calN(\bmu), \Delta_{n-1}(\bmu)): \bmu \in \Xi_{\rm train}\}$.
\vspace{0.05in}
\State 
\begin{minipage}{0.03\textwidth}
\rotatebox{90}{Multi-step greedy}
\end{minipage} 
\begin{minipage}{0.03\textwidth}
\rotatebox{90}{scan on $\Xi_{\rm Sur}$}
\end{minipage} 
\hspace{4mm}
\framebox[0.8 \textwidth]
{\begin{minipage}{0.8 \textwidth}
\While {$(\varepsilon > \varepsilon_{\mathrm{tol}})$ and $(\varepsilon > E_{\ell} \, \frac{1}{K_{\rm damp} \times (\ell + 1)})$}
\vspace{0.05in}
\For{each $\boldsymbol{\mu} \in \Xi_{\rm Sur}$}
\State Obtain RBM solution $u^{\mathcal{N}}_{n}(\boldsymbol{\mu}) \in X^\calN_n$ and error estimate ${\Delta_{n}}(\boldsymbol{\mu})$
\EndFor
\vspace{0.05in}
\State $\bmu^{n+1} = \underset{\bm{\mu} \in \Xi_{\rm Sur}}{\argmax} \Delta_{n}(\bm{\mu})$, $\varepsilon = \Delta_{n}(\bmu^{n+1})$.
\State Augment RB space $X^\calN_{n+1} = X^\calN_n \oplus \{u^\calN(\bmu^{n+1})\}$
\State Set $n \leftarrow n+1$, $N_\ell \gets N_\ell + 1$
\EndWhile
\end{minipage}
}
\vspace{0.05in}
\EndWhile
\end{algorithmic}
\caption{The Offline-enhanced Reduced Basis Method template. Algorithms for constructing the SPD $\Xi_{\rm Sur}$ are described in Sections \ref{sec:smm} and \ref{sec:cdm}.}\label{alg:oerbm}
\end{algorithm}

The other main contribution of our paper is the creation of two strategies for constructing the surrogate parameter domain $\Xi_{\rm Sur}$, which is the topic of the next two subsections. Our two procedures are the Successive Maximization Method (SMM) and the Cholesky Decomposition Method (CDM). Once they are described, we may use them in the algorithmic template that fully describes Offline-enhanced RB methods. SMM and CDM are intrinsically different in their construction, yet our numerical experiments show that they both work very well, accelerating the offline portion of the RBM algorithm significantly without sacrificing accuracy for the examples we have tested.

We make three remarks concerning the Algorithm template:
\begin{itemize}[itemsep=0pt]
  \item {\bf Motiviation for constructing the SPD} In a standard sweep of \eqref{eq:greedy-selection} to identify $\bmu^\ast$ from $\Xi_{\rm train}$ that maximizes the error estimate $\Delta$, we actually must compute $\Delta(\bmu)$ for \textit{all} $\bmu \in \mathcal{D}$. Standard RB algorithms discard this information upon identifying $\bmu^\ast$. However, this is valuable, quantitative information about $\left\|e_N(\bm{\mu})\right\| = \left\|u^{\mathcal{N}}(\bm{\mu}) - u_{N}^{\mathcal{N}}(\bm{\mu})\right\|$ for {\em any} $\bmu \in \mathcal{D}$. 
     Construction of $\Xi_{\rm Sur}$ attempts to utilize this information that was otherwise discarded to identify not just $\bmu^\ast$, but a collection of parameters that can describe the landscape of $\Delta(\cdot)$. In other words, we gauge the accuracy of the reduced solution in $X^{\mathcal{N}}_{n+1}$ for all $\mu$, and trim from $\Xi_{\rm train}$ those parameters whose corresponding solutions are deemed good enough. Roughly speaking, we set
    \begin{equation}\label{eq:Sur-qualitative}
      \Xi_{\rm Sur} := \left\{ \bmu: u_{n+1}^\calN(\bmu) \mbox{ is predicted to be ``inaccurate''}\right\}.
    \end{equation} 
    Note that one can mathematically devise adversarial scenarios where such a procedure can discard values in $\Xi_{\rm train}$ that later will be important. However, 
    the outer loop of the template is designed so that we reconsider any parameter values that may have been discarded at one point. The goal is to construct $\Xi_{\rm Sur}$ in a balanced way: A strict definition of ``inaccurate" in \eqref{eq:Sur-qualitative} makes $\Xi_{\rm Sur}$ too large and no savings is gained; a lax definition chooses too few values for $\Xi_{\rm Sur}$ and the RB surrogate will not be accurate. 
  \item {\bf Stopping criteria for the SPD} On outer loop round $\ell$, we repeatedly sweep the current SPD $\Xi_{\rm Sur}$ after it is constructed until 
    \begin{equation*}
      \max_{\bmu \in \Xi_{\rm Sur}} \Delta_k(\bmu) \le E_{\ell}\, \frac{1}{((\ell+1) \times K_{\rm damp})},
    \end{equation*}
    where $E_{\ell}$ is the starting (global) maximum error estimate for this outer loop iteration. The damping ratio $\frac{1}{((\ell+1) \times K_{\rm damp})}$, enforces that the maximum error estimate over the SPD decreases by a controllable factor $K_{\rm damp}$; in this paper we take $K_{\rm damp}$ to be constant in $\ell$.  However, this damping ratio should be determined by the practitioner and the problem at hand. Taking $K_{\rm damp}$ as a constant works well in our test problems. 
  \item {\bf Cost of constructing $\Xi_{\rm Sur}$} The cost of constructing $\Xi_{\rm Sur}$ is an overhead cost for each outer loop of Algorithm \ref{alg:oerbm}. Therefore, we must formulate this construction so that the overhead cost is worth the effort. For example, if the cost of evaluating $\Delta(\cdot)$ at one value is $C$, and we select $N_\ell$ snapshots from $\Xi_{\rm Sur}$ at outer iteration $\ell$, then we attain cost savings when
    \begin{align*}
      \frac{\textrm{Cost of constructing $\Xi_{\rm Sur}$}}{C N_\ell} < |\Xi_{\rm train}| - |\Xi_{\rm Sur}|.
    \end{align*}
    This yields qualitative information about the efficiency of the method: when the cost of constructing $\Xi_{\rm Sur}$ is negligible, we may take a large $\Xi_{\rm Sur}$, but when this cost is large, we require a significant size reduction in order to amortize the initial investment.
\end{itemize}

\subsection{Successive Maximization Method}\label{sec:smm}

Our first approach for constructing the surrogate parameter domain is the Successive Maximization Method (SMM). This procedure is motivated by the notion that the difference between the norm of the errors $\left | \lVert e(\bmu_1) \rVert_{X} - \lVert e(\bmu_2) \rVert_{X} \right |$ is partially indicative of the difference between the solutions. Computation of the true error norms is impractical, so like standard offline RBM procedures we leverage the {\em a posteriori} error estimate $\Delta_{N}(\bm{\mu})$ defined in \eqref{eq:error_estimator_en}. 

Suppose we have already selected $k$ snapshots; when selecting parameter value $k+1$ via \eqref{eq:greedy-selection}, we must compile the values $\Delta_k\left(\mathcal{D}\right) = \left\{ \Delta_k(\bmu) \;|\; \mu \in \mathcal{D}\right\}$. We use this collection to identify the surrogate parameter domain. From our argument that the values $\left\| e(\bmu)\right\|$ give us some indication about the actual solution, we equidistantly sample values from $\Delta_k(\mathcal{D})$ to construct the parameter domain. 

With $\epsilon_{\rm tol}$ the stopping tolerance for the RB sweep, let $\Delta_k^{\rm max} = \underset{{\bmu \in \Xi_{\rm train}}}{\max} \Delta_k(\bmu)$. We define $I_k^M$ as an {equi-spaced} set between $\epsilon_{\rm tol}$ and $\Delta_k^{\rm max}$:
\begin{align*}
  I_k^M = \left\{ {\nu_{k,m} \coloneqq \,} \epsilon_{\rm tol} + (\Delta_k^{\rm max} - \epsilon_{\rm tol}) \frac{m}{M}: m = {0}, \dots, {M - 1}\right\}.
\end{align*}
Roughly speaking, we attempt to construct $\Xi_{\rm Sur}$ as $\Xi_{\rm Sur} = \Delta_k^{-1}\left(I_k^M\right)$. Rigorously, we use
\begin{align*}
\Xi_{\rm Sur} &= \left\{ \bmu_{k,m}: \bmu_{k,m} = \mathrm{argmin}_{\bmu \in \Xi_{\rm train}} \left\{ \Delta_k(\bmu) - \nu_{k,m} \hskip 5pt \textrm{such that} \hskip 5pt \Delta_k(\bmu) \geq \nu_{k,m} \right\} \right\}.
\end{align*}
Note that we have $|\Xi_{\rm Sur}| {{\le}} M$ by this construction. 

\subsection{Cholesky Decomposition Method}\label{sec:cdm}
For the second approach, we consider the (scaled) Gramian matrix $G$ comprised of pairwise inner products of error vectors $e_k(\bmu)${.} I.e.,
\begin{align*}
  G_{i,j} &= \frac{(e(\bmu_i), e({\bmu_j}))_{X}}{\lVert e(\bmu_i) \rVert_{X} \lVert e({\bmu_j})\rVert_{X}}, & \bmu_i, \bmu_j &\in \Xi_{\rm train}{.}
\end{align*}
The matrix $G$ is positive semi-definite, and thus admits a (pivoted) Cholesky decomposition. We suppress notation indicating that $G$ depends on the current number of snapshots $k$.

Our approach here is to apply the {\em pivoted Cholesky decomposition} \cite{MR2899254} of the matrix $G$. This decomposition of $G$ orders the elements of $\Xi_{\rm train}$ according {to} the pivots. We identify the surrogate domain $\Xi_{\rm Sur}$ as the first $M$ pivots (parameter values) selected by this procedure. 

Since obtaining the error vectors $e(\bmu)$ is as expensive as solving for the truth approximation, we have to approximate these vectors. A linear algebraic way to write the Galerkin system \eqref{eq:update_problem} is 
\begin{align*}
  \mathbb{A}(\bmu) \bm{u}^{\calN}(\bmu) &= \bm{f}(\bmu), & \mathbb{A} &\in \R^{\calN \times \calN}
\end{align*}
where $\mathbb{A}$, $\bm{u}^{\calN}$, and $\bm{f}$ are discretization vectors associated to $a(\cdot, \cdot; \bmu)$, $u^{\calN}$, and $f(\cdot;\bmu)$, respectively. With this notation, we have
\[
  \bm{e}(\bm{\mu}) = \bm{u}^{\mathcal{N}}(\bm{\mu}) - \bm{u}^{\mathcal{N}}_{N}(\bm{\mu}) = \mathbb{A}^{-1}_{\mathcal{N}}(\bm{\mu}) \bm{r}(u_N^\calN(\bmu); \bmu), 
\]
where the residual vector $\bm{r}(\cdot; \bmu)$ is defined as $\bm{r}(\bm{v}; \bmu) = \bm{f} - \mathbb{A}_{\mathcal{N}}(\bm{\mu}) \bm{v}$. We propose to approximate the unknown $\mathbb{A}^{-1}_{\mathcal{N}}(\bm{\mu})$ by  
\begin{equation}\label{eq:A_inverse_approx}
 \widetilde{\mathbb{A}}{_\calN}^{-1} (\bm{\mu}) \coloneqq \sum_{m = 1}^{Q}u^{\mathcal{N}}_{Qm}(\bm{\mu}) \mathbb{A}{_\calN}^{-1}(\bm{\mu}^{m}),
\end{equation} 
where  $\left\{u^{\mathcal{N}}_{Qm}(\bm{\mu})\right\}_{m = 1}^{Q}$ are the RB coefficients for $u_Q^\calN(\bmu)$ defined in \eqref{eq:rbsolution}. Since $ ({\mathbb{A}}{_\calN}^{-1}(\bm{\mu}) -  \widetilde{\mathbb{A}}{_\calN}^{-1}(\bm{\mu})) f = u^\calN(\bmu) - u_N^\calN (\bmu)$, we argue that this approximation is reasonable. 

The approximation of $e(\bm{\mu})$ can be expressed as
$\widetilde{e}(\bmu) = \bigg(\sum_{m = 1}^{Q}u^{\mathcal{Q}}_{Nm}(\bm{\mu}) \mathbb{A}^{-1}_{\mathcal{N}}(\bm{\mu}^{m})\bigg)r(u_{N}^{\mathcal{N}}(\bm{\mu}))$ which admits an affine decomposition, and the Gramian matrix $G$ is approximated by
\begin{equation}\label{eq:error_matrix}
\widetilde{G}_{ij} = (\widetilde{e}(\bm{\mu}_{i}), \widetilde{e}(\bm{\mu}_{j}))/(||\widetilde{e}(\bm{\mu}_{i})||_{X} \times ||\widetilde{e}(\bm{\mu}_{j})||_{X}),
\end{equation}
where $1 \le i, \, j \le N_{\rm train}$.

Note that we take $Q < N$ since the matrices $\mathbb{A}$ are of size $\calN$ and are thus performing algebraic manipulations on them is expensive when $N$ is large.

\subsection{Complexity analysis}\label{sec:oerbm-complexity}

In Appendix \ref{sec:app-offline-online} we see that the computational complexity for the offline portion of the classical algorithm has order
\begin{align*}
  \overbrace{\calN^2 N^2 Q_a}^{\mbox{Reduced solve preparation}} + \overbrace{N_{\rm train} W_\alpha + {\boldsymbol {N_{\rm train} (Q_a^2 N^3+ N^4)}} + N W_s}^{\mbox{Greedy sweeping}} + \overbrace{Q_a^2 N^3 W_m}^{\mbox{Estimator preparation}}.
\end{align*}
This cost is dominated by the boldface term in the middle, especially when $\Xi_{\rm train}$ is large. We denote this cost $\calC_{\rm orig}$,
\begin{align*}
  \calC_{\rm orig} \coloneqq N_{\rm train} (Q_a^2 N^3 + N^4).
\end{align*}
This is the portion of the offline cost that our Offline-Enhanced RBM is {aiming} to reduce.

Suppose we have computed $n$ snapshots. Then the cost for assembling and solving the RB system for one given parameter value $\bmu$ is of order $n^2 Q_a + n^3$, while the cost for calculating the error certificate is of order $n^2 Q_a^2$. Therefore the total cost for one instance, denoted by $c(n)$, is of order $n^2 Q_a^2 + n^3$. This means that the complexity for the classical RBM to sweep over $N_{\rm train} \coloneqq |\Xi_{\rm train}|$ parameter values is
\begin{align*}
N_{\rm train} \sum_{n = 1}^{N} \, c(n) = \calO\left( N_{\rm train} (Q_a^2 N^3 + N^4) \right).
\end{align*}
To better analyze the cost of our Offline-enhanced approaches, we denote the {\em cumulative} number of chosen parameter values after the $j$-th outer loop iteration by 
\[
t_j = 1 + \sum_{k=1}^j \, N_k 
\]
We note that $t_0 = 1$ because in standard RB algorithms the first parameter value is randomly chosen before starting the greedy algorithm. If the RB procedure given by Algorithm \ref{alg:oerbm} terminates after $\ell$ outer loop iterations with a total of $N$ snapshots, then we have $t_{\ell} = N$. 
The cost of the Offline-enhanced approaches corresponding to the dominating cost of the classical approach $\calC_{\rm orig}$ is
\begin{align*}
\overbrace{\calC_{\rm oe}}^{\mbox{Total cost}} & = \overbrace{N_{\rm train}\,\sum_{j=0}^{\ell}\, c(t_j)}^{\mbox{ One-step greedy scans on $\Xi_{\rm train}$}} + 
  \overbrace{\sum_{j=1}^{\ell} \, M_j \left( \sum_{n = t_{j-1}+1}^{t_j- 1} c(n) \right)}^{\mbox{Multi-step greedy scans on $\Xi_{\rm Sur}$}}\\
& \lesssim N_{\rm train} \ell c(N) + M_{\rm max} \sum_{n=1}^{N} c(n), \quad {\mbox{ with }} M_{\rm max } \coloneqq \max_{j = 1, \dots, \ell} M_j\\
& \lesssim N_{\rm train} \ell (Q_a^2 N^2 + N^3) + M_{\rm max} (Q_a^2 N^3 + N^4).
\end{align*}
We conclude that the dominating parameter sweeping costs between the classical and offline-enhanced approaches satisfy
\begin{equation}
\calC_{\rm oe} < \calC_{\rm orig} \left( \frac{\ell}{N} + \frac{M_{\rm max}}{N_{\rm train}} \right).
\label{eq:costratio}
\end{equation}
We make some remarks concerning this cost analysis:
\begin{itemize}[itemsep=0pt]
  \item {\bf Potential savings} --- Since $\frac{M_{\rm max}}{N_{\rm train}}$ is negligible especially for the cases of our concern when the parameter dimension is high, \eqref{eq:costratio} demonstrates that the savings is roughly $\frac{\ell}{N}$. 
  \item {\bf Surrogate domain construction for SMM-RBM} --- The additional cost for the Offline-enhanced approaches is the construction of the surrogate parameter domain. For SMM-RBM, this cost is essentially negligible. This SMM surrogate domain construction cost is mainly dependent on the cost of evaluating the error certificate $\Delta_k$, but this cost has already been cataloged in the analysis above. In practice, the surrogate domain construction amounts to a quick sorting of these certificates which results in a cost of $\calO(N_{\rm train}\log (N_{\rm train}))$. While this cost does depend on $N_{\rm train}$, it is much smaller than any of the terms in, e.g., $\calC_{\rm orig}$.
  \item {\bf Surrogate domain construction for CD-RBM} --- CD-RBM entails a sequence of (pivoted) Cholesky Decomposition steps applied to the approximate error Gramian matrix $\widetilde{G}$, \eqref{eq:error_matrix}. For the decomposition algorithm, it suffices to just supply the approximate errors $\widetilde{e}(\bmu)$ without constructing the full matrix $\widetilde{G}$. Therefore, the cost is primarily devoted to computing these approximate errors for all $\bmu$.  Evaluating these error functions can be accomplished in an offline-online way detailed in Appendix \ref{sec:app-offline-online-CDG}. We summarize here the total cost for the SPD construction.
\begin{align*}
\overbrace{\mathcal{N}^2QN_{\rm RB}Q_{a}}^{\mbox{ Offline Preparation }} + 
\overbrace{N_{\rm train} \mathcal{N}QN_{\rm RB}Q_{a} \ell}^{\mbox{Approximate error calculation}} +
  \overbrace{\mathcal{N} \sum_{j=1}^{\ell} \left(n^j_{\rm cd}\right)^2}^{\ell \mbox{ runs of Pivoted CD algorithm}}
\end{align*}
where $n^j_{\rm cd}$ is the number of steps of the pivoted Cholesky decomposition for the $j$-th iteration. 
We see that this algorithm can be more costly than SMM-RBM because of the factor $\mathcal{N}QN_{\rm RB}Q_{a}$. However, we observe that it is still much faster than the classical version since this factor is notably smaller than $Q_a^2 N^3_{\rm RB} + N^4_{\rm RB}$. This is confirmed by our numerical examples presented in the next section.
\end{itemize}

\section{Offline-enhanced RBM: Numerical results}\label{sec:results}
In this section, we present numerical examples to illustrate the accuracy and efficiency enhancement of the proposed approaches compared to the conventional reduced basis method. 

\subsection{Test problems} We test the two algorithms, SMM-RBM and CD-RBM, on two standard diffusion-type problems. They vary substantially in terms of parameter dimension and the truth solver. Our results show that our offline-enhanced procedures work well in both of these cases, and suggest that the offline-enhanced strategies may be beneficial for a wider class of parameter spaces and truth solvers. 

\noindent{\bf Diffusion problem with two-dimensional parameter domain:}
\begin{equation}\label{eq:diff_prob}
(1+\mu_{1}x)u_{xx} + (1+\mu_{2}y)u_{yy} = e^{4xy} \quad {\rm on} \quad \Omega.
\end{equation}
Here $\Omega = [-1 ,1] \times [-1,1]$ and we impose homogeneous Dirichlet boundary conditions on $\partial \Omega$. The truth approximation is a spectral Chebyshev collocation method based on $\mathcal{N}_{x} = 35$ degrees of freedom in each direction, with $\mathcal{N}^2_{x} = \mathcal{N}$.   The parameter domain $\mathcal{D}$ for $(\mu_{1}, \mu_{2})$  is taken to be $[-0.99, 0.99]^2$. For the training set $\Xi_{\rm train}$ we discretize $\calD$ using a tensorial $160 \times 160$ Cartesian grid with 160 {equi-spaced} points in each dimension. \\

\noindent {\bf Thermal Block problem with nine-dimensional parameter domain:}
\begin{equation}\label{eq:diff_prob_2}
\begin{cases}
- \nabla . (a(x, \bmu) \nabla u(x, \bmu) ) = f \quad \text{on}  \quad \Omega,\\
u(x,\bmu) = g_{D} \quad \text{on} \quad \Gamma_{D}, \\
\frac{\partial u}{\partial n} = g_{N} \quad \text{on} \quad \Gamma_{N}.
\end{cases}
\end{equation}
Here $\Omega = [0, 1] \times [0, 1]$ which is partitioned into $9$ blocks $\cup_{i=1}^9 B_i = \Omega$, $\Gamma_{D}$ is the top boundary, and $\Gamma_{N}  = \partial \Omega   \setminus \Gamma_{D}$. The parameters $\mu_{i}, ~ 1 \le i \le 9$ denote the heat conductivities:
\begin{center}
\begin{tabular}{|c|c|c|}
\multicolumn{3}{c}{$\Gamma_{D}$} \\
\hline
$\mu_{7} (B_{7})$ & $\mu_{8} (B_{8})$ & $\mu_{9} (B_{9})$ \\
\hline
$\mu_{4} (B_{4})$ & $\mu_{5} (B_{5})$ & $\mu_{6} (B_{6})$ \\
\hline
$\mu_{1} (B_{1})$ & $\mu_{2} (B_{2})$& $\mu_{3} (B_{3})$ \\
\hline
\multicolumn{3}{c}{$\Gamma_{base}$} 
\end{tabular}
\end{center}
The diffusion coefficient $a(x, \bmu) = \mu_{i} ~\text{if} ~x \in B_{i}$. The parameter vector is thus given by $\bmu = (\mu_{1}, \mu_{2}, \dots, \mu_{9})$  in $\mathcal{D} = [0.1, 10]^9$. We take as the right hand side $f = 0, g_{D} = 0$,  $g_{N} = 1$ on the bottom boundary $\Gamma_{\rm base}$ and $g_{N} = 0$ otherwise. The output of interest is defined as the integral of the solution over $\Gamma_{\rm base}$
\begin{equation}\label{eq: output}
s(\bmu) = \int_{\Gamma_{\rm base}} u(x,\bmu)dx
\end{equation}
The truth approximation is obtained by FEM with $\mathcal{N} = 361$. A sufficient number of $N_{{\rm train}} = {20,000}$ samples are taken from randomly sampling within the parameter domain $\mathcal{D}$. The classical RB solver for this problem is provided by, and our Offline-enhanced approach is compared against, the RBmatlab package \cite{RBmatlab, RBmatlabPaper}\footnote{Available for download at {\texttt http://www.ians.uni-stuttgart.de/MoRePaS/software/index.html}}.\\

\subsection{Results}
We investigate the performance of the two offline-enhanced RB algorithms on the two test problems. The tuning parameters for Algorithm \ref{alg:oerbm} for both examples are shown in Table \ref{tab:para}.
\begin{table}
\renewcommand{\tabcolsep}{0.4cm}
\renewcommand{\arraystretch}{1.5}
\begin{center}
\begin{tabular}{c|c|c}
\hline
Parameter & SMM-RBM & CD-RBM\\
\hline
$M_\ell$ & $2 \times (\ell + 1)$ & $20 \times (\ell + 1) $\\
\hline
$K_{\rm damp}$ & 1 & 10\\
\hline
\end{tabular}
\end{center}
\caption{Offline-enhanced RBM parameters for the numerical results.}
\label{tab:para}
\end{table}

\begin{figure}
\begin{center}
\resizebox{\textwidth}{!}{
  \includegraphics[width=0.33\textwidth]{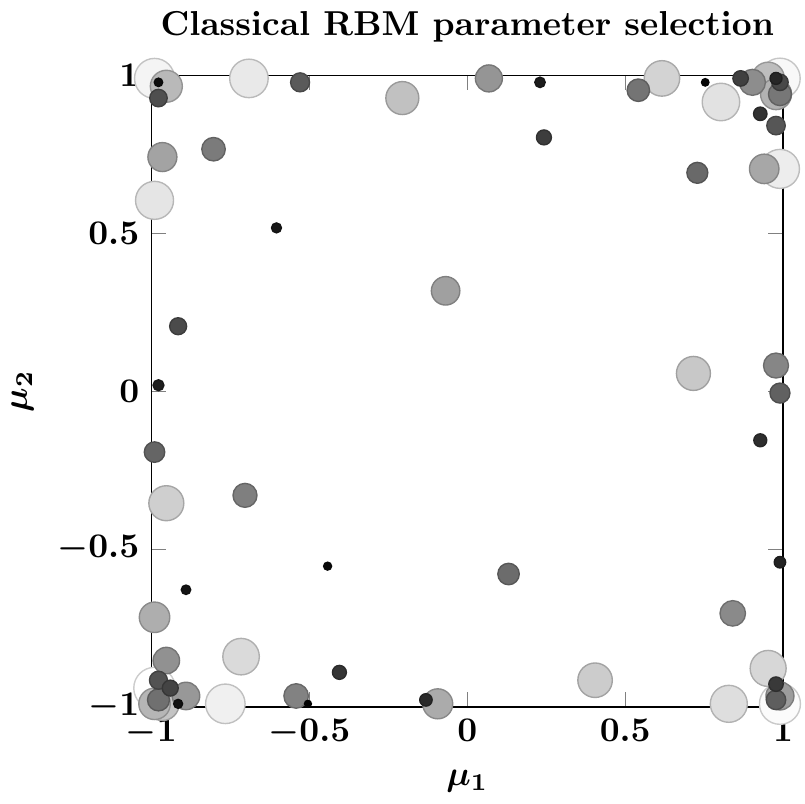}
  \includegraphics[width=0.33\textwidth]{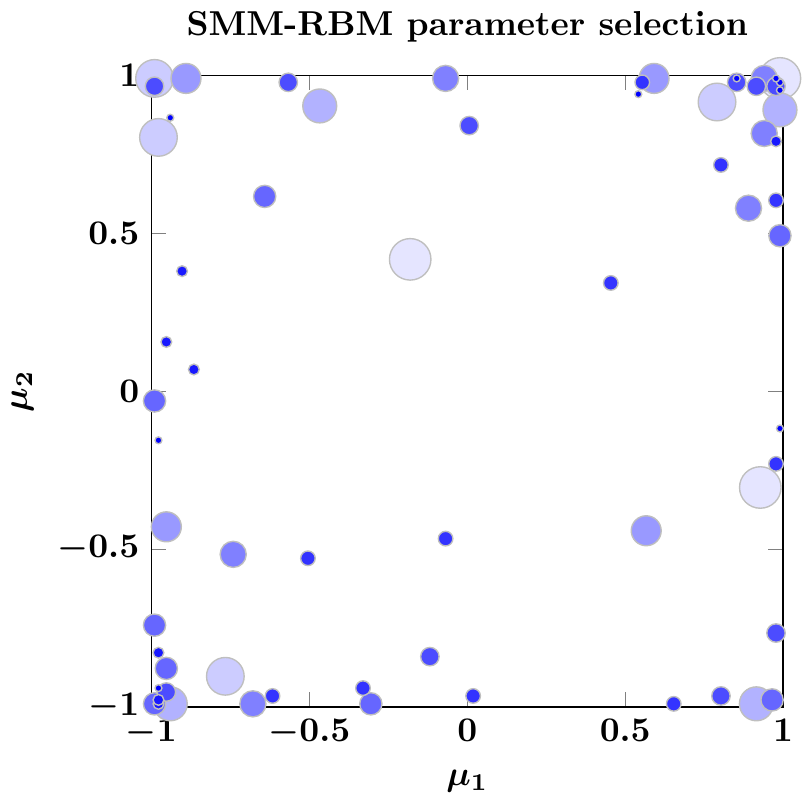}
  \includegraphics[width=0.33\textwidth]{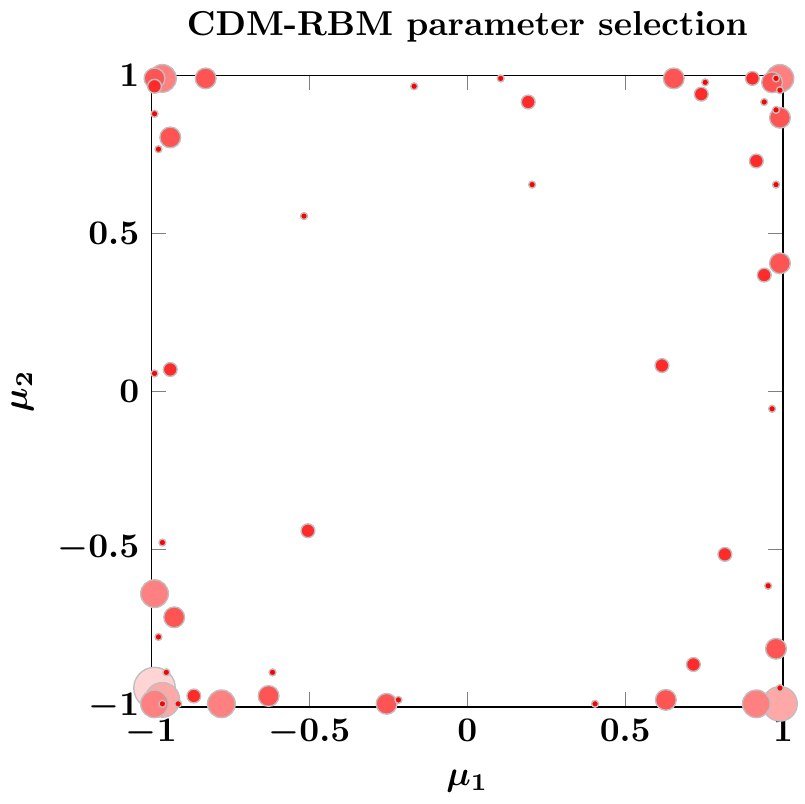}
}
\end{center}
\caption{Selected parameter values by a classical offline RBM algorithm (left), by the SMM-RBM algorithm (center), and by the CDM-RBM algorithm (right). Points with larger radius are chosen earlier in the sequence, and for the SMM and CDM plots points earlier in the sequence have greater transparency. For the SMM and CDM plots, all chosen parameters {within} a batch (an outer loop in Algorithm \ref{alg:oerbm}) have the same radius.}
\label{fig: RBM_scatter_prob1}
\end{figure}

\begin{figure}
\begin{center}
  \resizebox{!}{0.24\textheight}{
    \includegraphics[width=\textwidth]{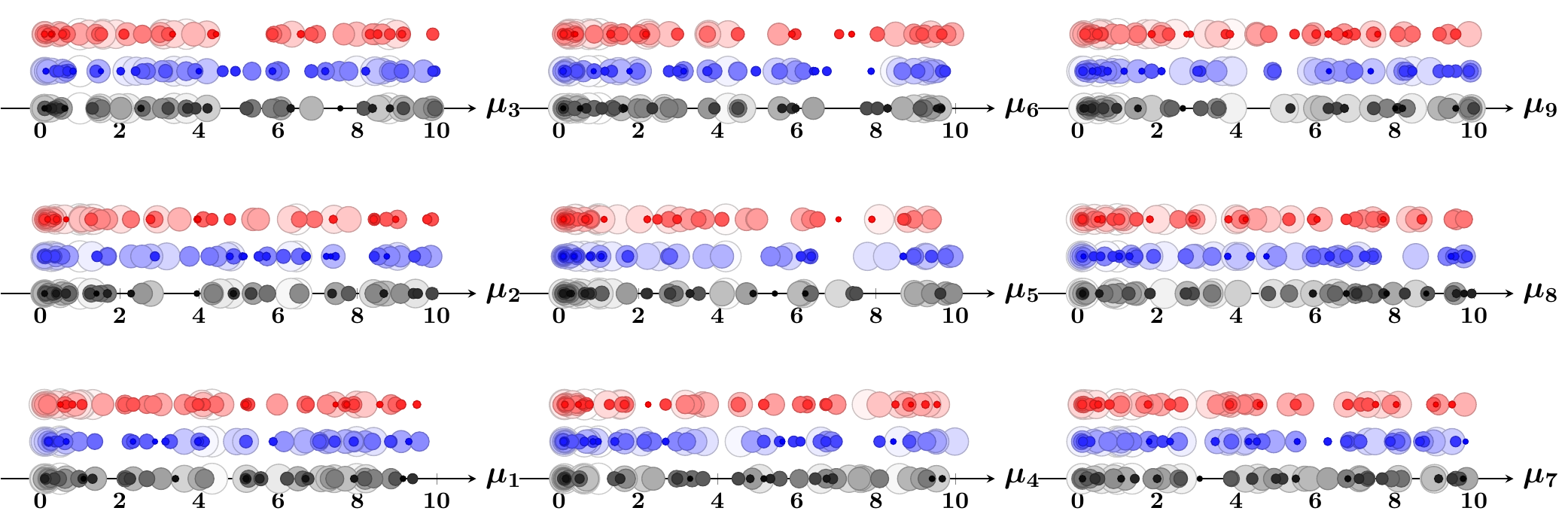}
  }
\end{center}
\caption{Parameter selections for the test problem 2 with a 9-dimensional parameter. One-dimensional scatter plots of the parameter selection are shown: classical RBM (black, bottom), SMM-RBM (blue, middle), and CDM-RBM (red, top). The size and transparency of the markers indicate the ordering of the sequence: points earlier in the sequence have a larger radius and are more transparent, points later in the sequence are smaller and more opaque.}
\label{fig:univariate-marginals-problem2}
\end{figure}

For the first test case with a two-dimensional parameter domain that is easy to visualize, we display the location 
of the selected parameter values in Figure \ref{fig: RBM_scatter_prob1}. On the {left} is that for the classical RBM. 
{The larger the marker, the earlier the parameter is picked.} {At the middle and on the right} are the parameter sets selected by SMM-RBM and CD-RBM respectively. For these two, the {more transparent} the marker, the earlier it is picked. 
A group of parameter values chosen at the same step {have the same radius}. 
Figure \ref{fig:univariate-marginals-problem2} displays the one-dimensional marginal scatter plots of the 9-dimensional point selection for the second test case.  
The two-dimensional marginal scatter plots are shown in the Appendix, 
in Figures \ref{fig:bivariate-marginals-problem2-rb}, \ref{fig:bivariate-marginals-problem2-smm}, and \ref{fig:bivariate-marginals-problem2-cdm} for the classical RB offline algorithm, the SMM-RBM algorithm, and the CDM-RBM algorithm, respectively.


\begin{figure}
\begin{center}
  \resizebox{\textwidth}{!}{\includegraphics{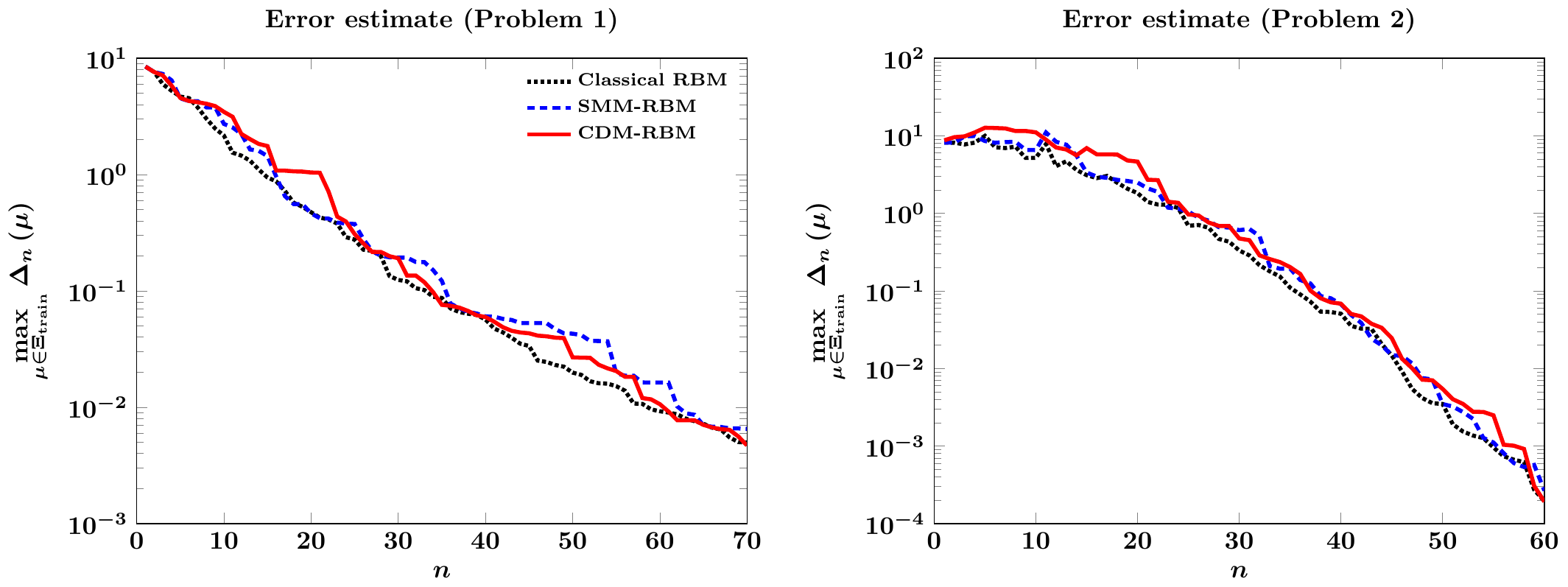}}
\end{center}
\begin{center}
  \resizebox{\textwidth}{!}{\includegraphics{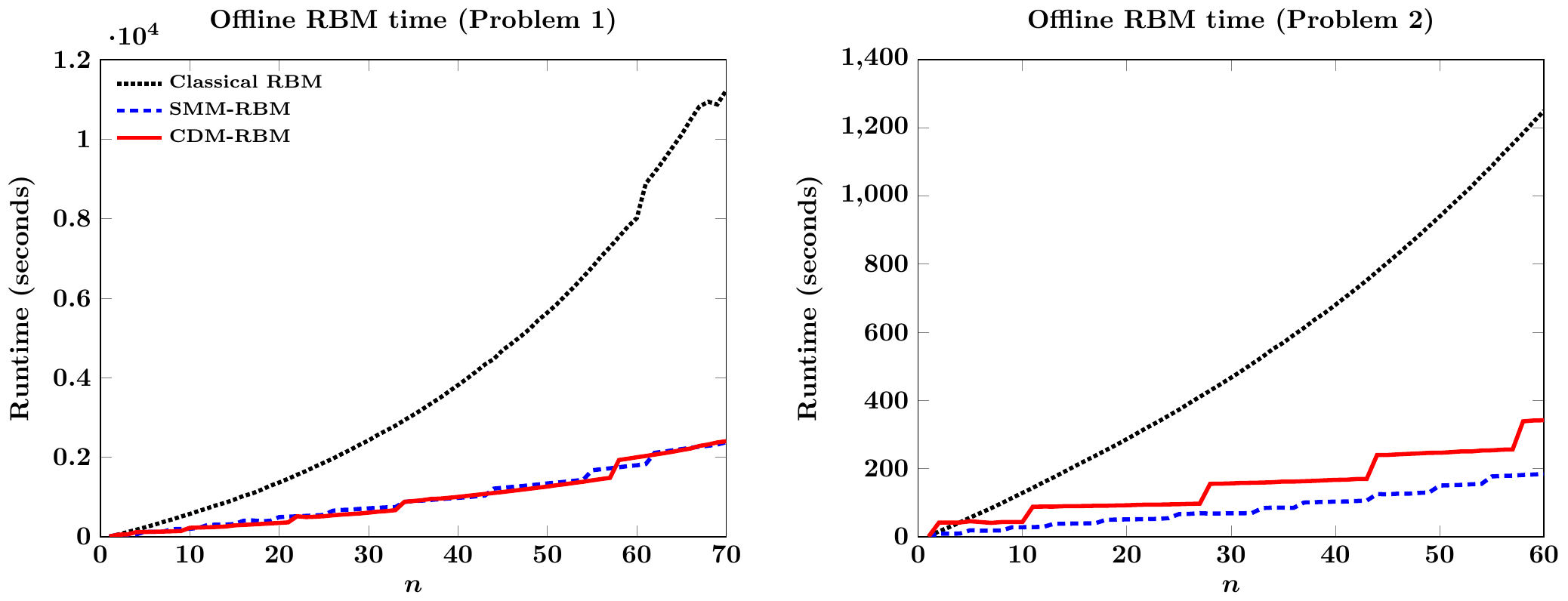}}
\end{center}
\begin{center}
  \resizebox{\textwidth}{!}{\includegraphics{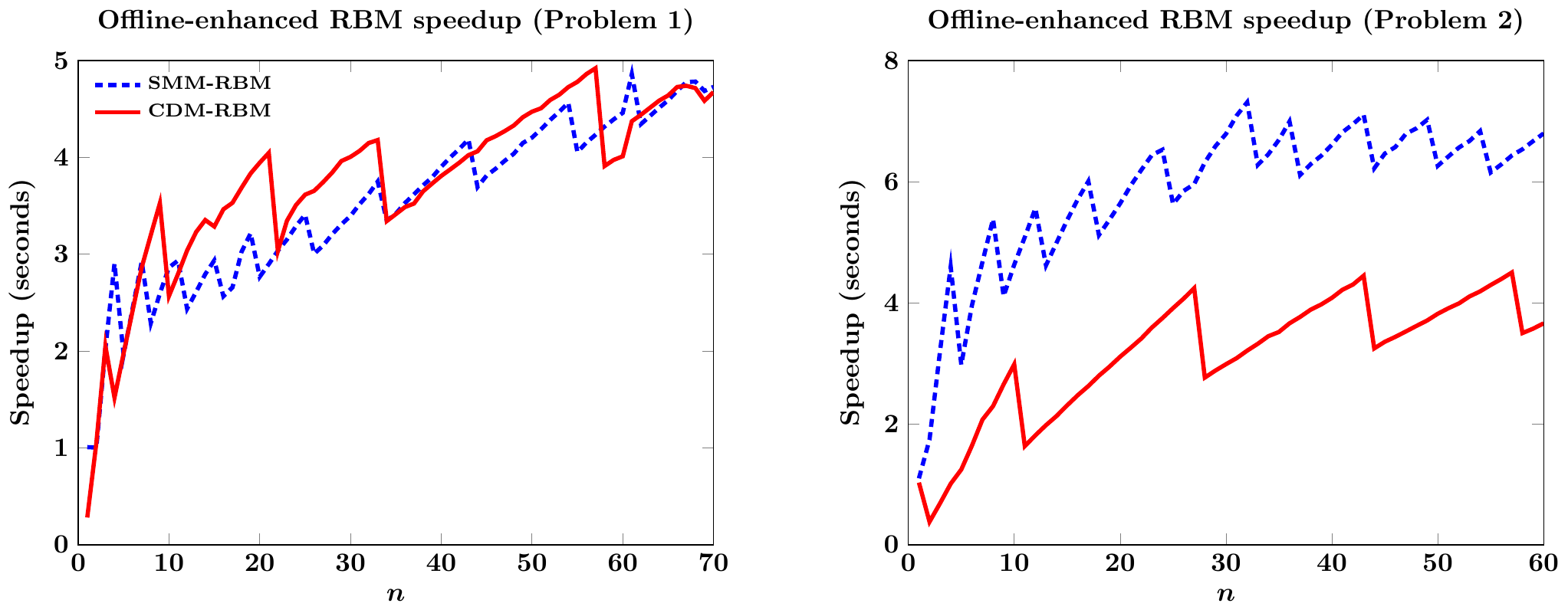}}
\end{center}
\caption{Convergence and speedup of offline-enhanced RBM algorithms as a function of the total number of snapshots $n$. Error estimate convergence (top), computational runtime (middle), and speedup factor (bottom). Left-hand plots correspond to the test problem 1, and right-hand plots for test problem 2.}
\label{fig: ET_prob1}
\end{figure}

The accuracy and efficiency of the new algorithms are shown in Figure \ref{fig: ET_prob1}. 
We see clearly that the {\em a posteriori} error estimate is 
converging exponentially for both SMM-RBM and CD-RBM, in the same fashion as the classical version of RBM. 
This shows that our accelerated algorithm does not appear to suffer accuracy degradation for these examples. In addition, we see a factor of $3$-to-$6$ times runtime speedup.

\begin{figure}
\begin{center}
  \includegraphics[width=0.49\textwidth]{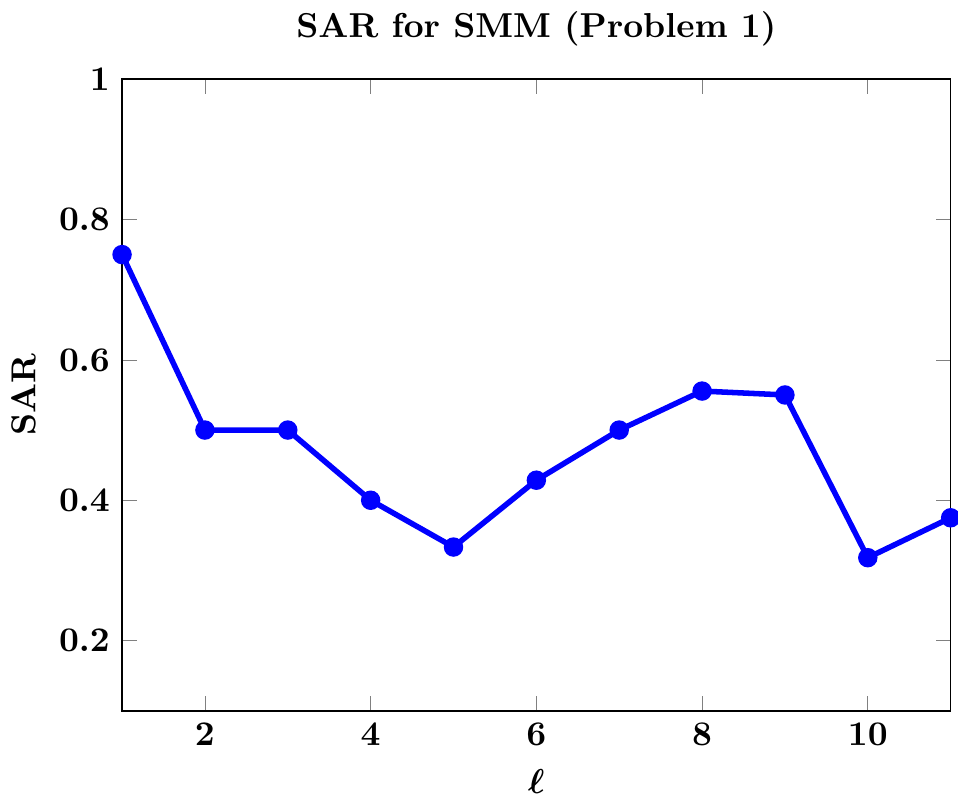}
  \includegraphics[width=0.49\textwidth]{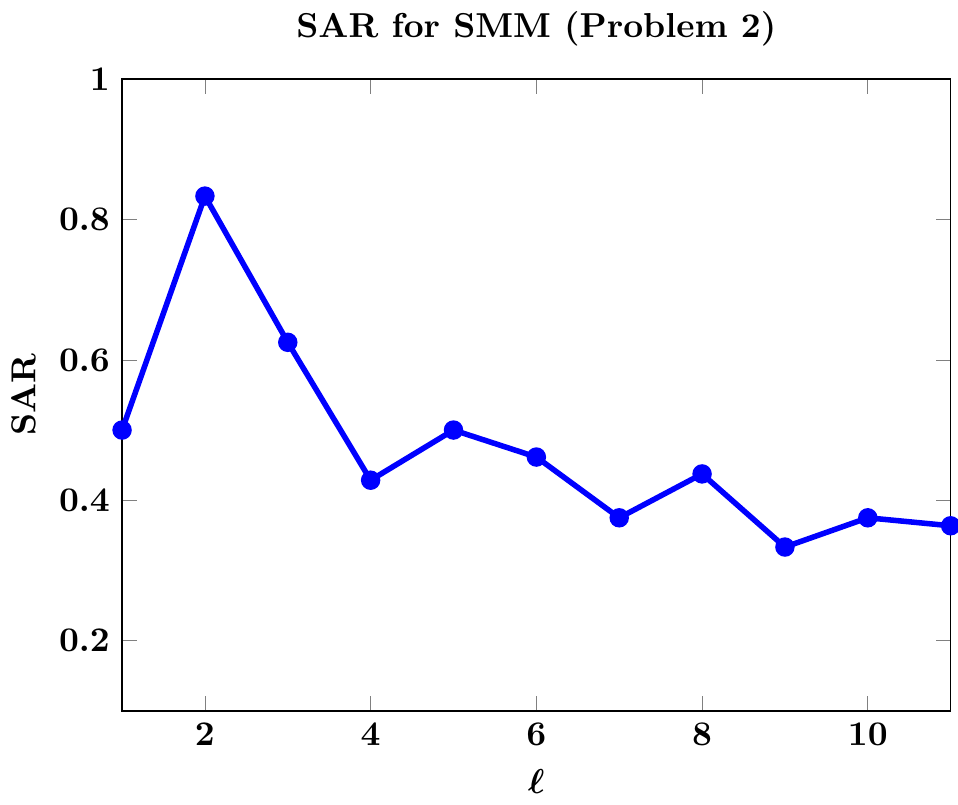}
\end{center}
\begin{center}
  \includegraphics[width=0.49\textwidth]{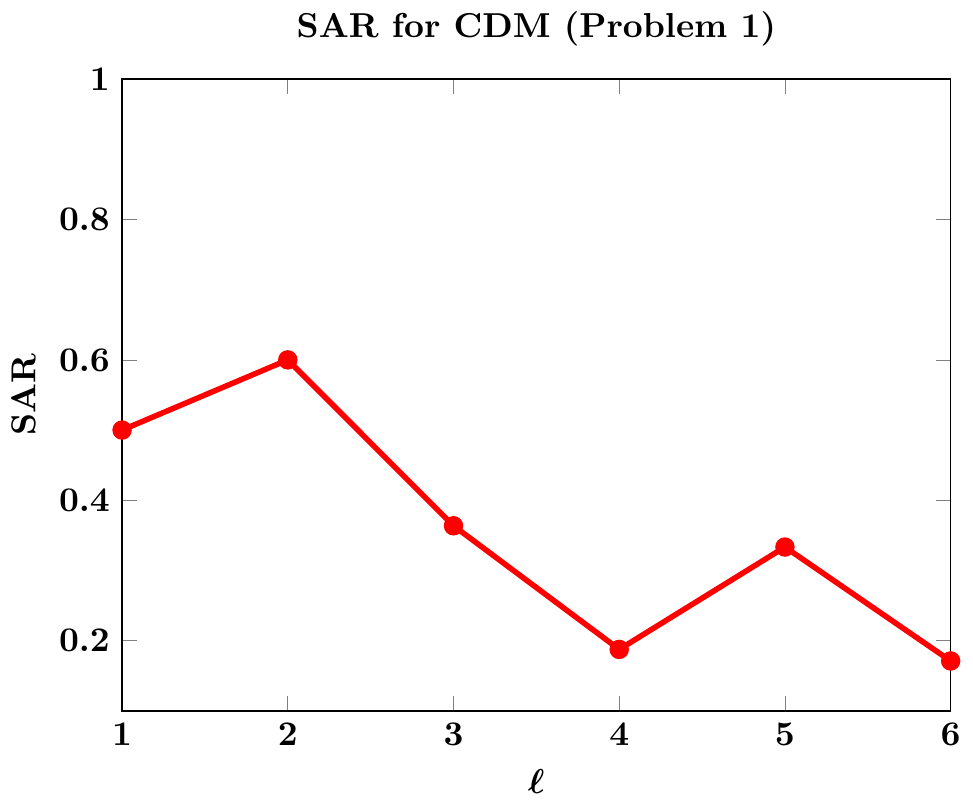}
  \includegraphics[width=0.49\textwidth]{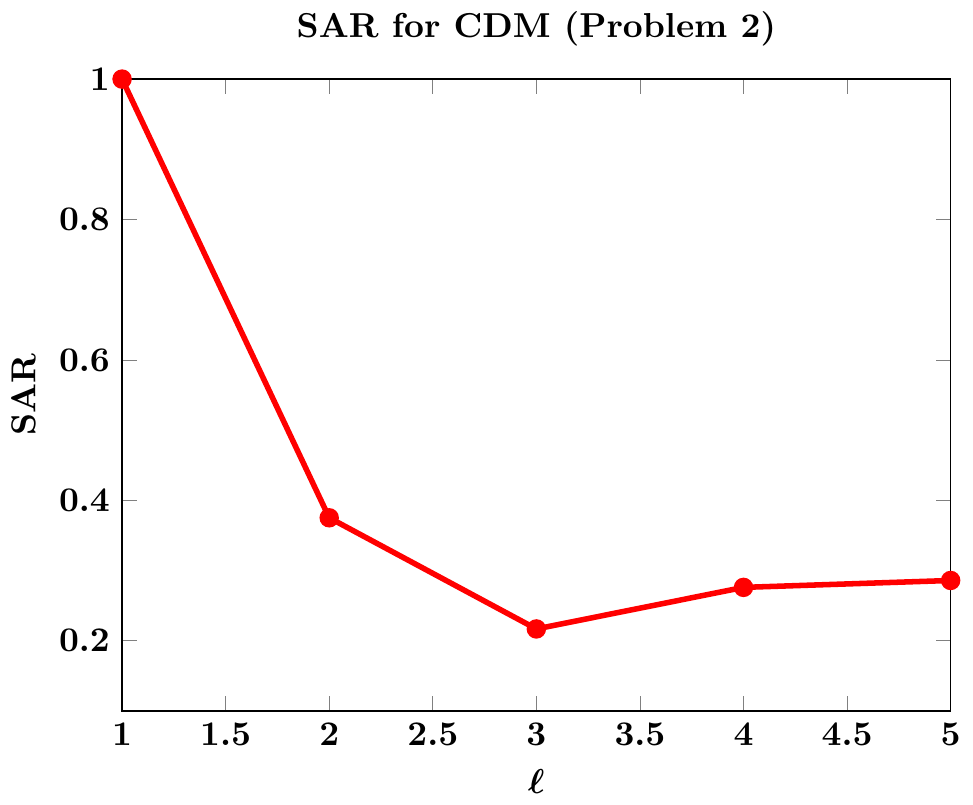}
\end{center}
\caption{The Surrogate Acceptance Ratio (SAR) for the offline-enhanced RBM methods, SMM (top) and CDM (bottom). Test problem 1 is shown on the left, and test problem 2 on the right.}
\label{fig:pickingratio}
\end{figure}

Finally, to reveal the effectivity of the construction of the surrogate parameter domains, we plot in Figure \ref{fig:pickingratio} the \textit{Surrogate Acceptance Ratio}
\begin{align*}
  \textrm{SAR} (\ell) \coloneqq \frac{N_\ell}{M_\ell}
\end{align*}
as a function of the outer loop iteration index $\ell$. This ratio quantifies how much of the surrogate domain is added to the snapshot parameter set at each outer loop iteration. Large ratios suggest that our construction of the surrogate domain effectively emulates the entire training set. In Figure \ref{fig:pickingratio} we see that a significant portion (on average approximately $40\%$ for the first case and $30\%$ for the second case) of the surrogate parameter domain is chosen by the greedy algorithm before continuing into another outer loop. 
We recall that the second case has a $9$-dimensional parameter domain, and so our offline-enhanced RBM procedure can effectively choose surrogate domains even when the parameter dimension is large. The relatively large values of the SAR result in the computational speedup observed in Figure~\ref{fig: ET_prob1}.

\section{Concluding Remarks}
We proposed an offline-enhanced reduced basis methods for building reduced-order models; RBM algorithms invest significant resources in an offline stage by studying a finite training set and judiciously choosing snapshots from this training set. Our novel approach substitutes the original training set with an adaptively constructed surrogate domain that is much smaller in size, and thus reduces the computational time expended in the offline portion of the RBM algorithm. (Our algorithm leaves the online portion of RBM algorithms unchanged.)

We provide two approaches to identify and construct surrogate domains using two different perspectives: the SMM-RBM strategy constructs a surrogate domain by uniformly sampling parameters on the range of the \textit{a posteriori} error estimate; the CD-RBM strategy analyzes the angle between two approximate error vectors at different locations in parameter space. Like RBM in general, our approaches are particularly useful in computing many-query reliable solutions parametrized PDE having a large number {of} random inputs. We have demonstrated the computational efficiency of our proposed methods compared against standard reduced basis method for two steady-state diffusion problems. The application of the offline-enhanced reduced basis method to more general problems with high dimensional parameter domains is ongoing research.

\appendix

\section{Classical RBM specifics: greedy algorithms, efficiency, and operational count}

\label{appendix}

This appendix contains the mathematical and algorithmic portions of RBM algorithms that are not directly the subject of this manuscript. These specifics are well-known in the RBM literature and community, and we include this appendix mainly for completeness of this manuscript. Section \ref{sec:app-greedy} discusses the mathematical justification for why the greedy procedure \eqref{eq:greedy-selection} is a good selection of parameter snapshots. Section \ref{sec:app-error-estimate} discusses efficient computation of the error estimate $\Delta_n(\bmu)$; the main goal is to compute this estimate with computational complexity that is independent of the truth solution complexity $\calN$. Section \ref{sec:app-offline-online} gives an overview of the RBM procedure, and quantifies the computational complexity of the RBM algorithm. Careful scrutiny of this operational count illustrates why RBM algorithms can simulate parameterized problems with $\calN$-independent complexity in the online phase of the algorithm.

Finally, Section \ref{sec:app-offline-online-CDG} discusses an efficient methodology to compute entries of the approximate Gramian $\widetilde{G}$ used by \eqref{eq:error_matrix} in the CDM algorithm. This procedure is a relatively straightforward application of the offline-online decomposition already employed by RBM algorithms.

\subsection{Greedy and Weak Greedy Algorithms}\label{sec:app-greedy}
The best $N$-dimensional RB space $X_{N}^{\mathcal{N}}$ in $X^{\cal}$ among all possible $N$-dimensional subspaces of the solution manifold $u\left(\cdot; \calD \right)$ is in theory the one with the smallest Kolmogorov $N$-width $d_N$ \cite{MR774404}:
\begin{align}\label{eq:kolmogorov-n-width}
  d_N \left[ u\left(\cdot; \calD\right) \right] \coloneqq \inf_{\substack{X_N \subset X^{\calN} u\left(\cdot; \Gamma\right) \\ \dim X_N = N}}\;\; \sup_{\mu \in \calD}\;\; \inf_{v \in X_N} \left\| u(\cdot, \mu) - v \right\|_X
\end{align}
The identification of an exact-infimizer for the outer ``inf'' is usually infeasible, but a prominent approach is to employ a greedy strategy 
which locates this $N$-dimensional space hierarchically. 
A first sample set $S_{1} = \{\bm{\mu}^{1}\}$ is identified by randomly selecting $\bmu^1$ from $\Xi_{\rm train}$; its associated reduced basis space $X_{1}^{\mathcal{N}} = \text{span}\{u^{\mathcal{N}}(\bm{\mu}^{1})\}$ is likewise computed. Subsequently parameter values are greedily chosen as sub-optimal solutions to an $L^{2}(\Xi_{\rm train}; X)$ optimization problem \cite{RozzaHuynhPatera2008}: for $N = 2, \dots, N_{\rm max}$, we find 
\begin{equation}
\label{eq:opt_problem}
\bm{\mu}^{N} = \underset{\bm{\mu} \in \Xi_{\rm train}} {\text{argmax} }~ || u^{\mathcal{N}}(\bm{\mu}) - u^{\mathcal{N}}_{N-1}(\bm{\mu})||_{X^{\calN}}
\end{equation}
where $u^{\mathcal{N}}_{N-1}(\bm{\mu})$ is the RB solution \eqref{eq:rbsolution} in the current $(N-1)$-dimensional subspace. 
Direct calculation of $u^\calN (\bmu)$ to solve this optimization problem over all $\bmu$ is impractical.  Therefore, a weak greedy algorithm is usually employed where we replace the error $\lVert u^{\mathcal{N}}(\bm{\mu}) - u^{\mathcal{N}}_{N-1}(\bm{\mu})\rVert_{X}$ by an inexpensive and computable \textit{a posteriori} bound $\Delta_{N-1}$ (see the next section).
After identifying $\bmu^N$, the parameter snapshot set and the reduced basis space are augmented, $S_{N} = S_{N-1} \cup \{\bm{\mu}^{N}\} ~\text{and} ~ X^{\mathcal{N}}_{N} = X^{\mathcal{N}}_{N-1} \oplus \{u(\bm{\mu}^{N})\} $, respectively.

\subsection{A posteriori error estimation}\label{sec:app-error-estimate}
The design of an effective {\em a posteriori} error bound $\Delta_{N-1}$ is crucial for the reliability of the reduced basis space constructed according to the weak greedy algorithm discussed in the previous section. Toward that end, we reconsider the numerical schemes for the truth approximation \eqref{eq:update_problem} and for the RB solution \eqref{eq:reduced_system}. Defining the error $e_N(\bm{\mu}) := u^{\mathcal{N}}(\bm{\mu}) - u^{\mathcal{N}}_{N}(\bm{\mu}) \in X^{\mathcal{N}}$, linearity of $a$ yields the following error equation:
\begin{equation}
\label{eq:residual_eq}
a(e_N(\bm{\mu}), v; \bm{\mu}) = r_N(v; \bm{\mu}) \quad \forall v \in X^\calN,
\end{equation}
with the residual $r_N(v; \bm{\mu}) \in (X^{\mathcal{N}})'$ (the dual of $X^{\calN}$) is defined as $f(v; \bm{\mu}) - a(u_{N}^{\mathcal{N}}(\bm{\mu}), v; \bm{\mu})$. The Riesz representation theorem and the Cauchy-Schwarz inequality implies that $\lVert e_N(\bm{\mu})\rVert_{X} \le \frac{\lVert r_N(\cdot; \bm{\mu})\rVert_{(X^{\mathcal{N}})'}}{\alpha_{LB}^{\mathcal{N}}(\bm{\mu})}$, 
where $\alpha^{\mathcal{N}}(\bm{\mu}) = \underset{w \in X^{\mathcal{N}}}{\inf} \frac{a(w, w, \bm{\mu})}{||w||^2_{X}}$ is the stability (coercivity) constant for the elliptic bilinear form $a$. This implies that we can define the {\em a posteriori} error estimator for the solution as
\begin{equation}\label{eq:error_estimator_en}
  \Delta_{N}(\bm{\mu}) = \frac{\lVert r_N(\cdot; \bm{\mu})\rVert_{(X^{\mathcal{N}})'}}{\alpha^{\mathcal{N}}_{LB}(\bm{\mu})} \geq \left\| e_N \left(\bmu\right) \right\|_{X^{\calN}}
\end{equation}
 The efficiency of computing the {\em a posteriori} error estimation relies on that of the lower bound of the coercivity constant $\alpha^{\mathcal{N}}_{LB}(\bm{\mu})$ as well as the value $\lVert r_N(\cdot; \bm{\mu})\rVert_{(X^{\mathcal{N}})'}$ for $\forall \bm{\mu} \in \mathcal{D}$. The coercivity constant $\alpha^{\calN}$ can be nontrivial to compute, but there are constructive algorithms to address this \cite{DGSA2007,HuynhSCM,Chen2015_NNSCM,HKCHP}. The residual is typically computed by the RBM offline-online decomposition, which is the topic of the next section. 

\subsection{Offline-Online decomposition}\label{sec:app-offline-online}

The last component of RBM that we plan to review in this section is the Offline-Online decomposition procedure \cite{RozzaHuynhPatera2008}. The complexity of the offline stage depends on $\calN$ which is performed only once in preparation for the subsequent online computation, whose complexity is independent of $\calN$. It is in the $\calN$-independent online stage where RBM achieves certifiable orders-of-magnitude speedup compared with other many-query approaches. The topic of this paper addresses acceleration of the offline portion of the RBM algorithm. In order to put this contribution of this paper in context, in this section we perform a detailed complexity analysis of the decomposition.

We let $N_{\rm train} = |\Xi_{\rm train}|$ denote the cardinality (size) of $\Xi_{\rm train}$; $N \le N_{\rm max}$ is the dimension of the reduced basis approximation computed in the offline stage. Computation of the the lower bound $\alpha_{LB}^{\mathcal{N}}(\bm{\mu})$ is accomplished via the Successive Constraint Method \cite{DGSA2007}. 

During the online stage and for any new $\bmu$, the online cost of evaluating $\alpha_{LB}^{\mathcal{N}}(\bm{\mu})$ is negligible, but we use $W_{\alpha}$ to denote the {\em average} cost for evaluating these values over the training set $\Xi_{\rm train}$ (this includes the offline cost). $W_{s}$ is the operational complexity of solving problem \eqref{eq:update_problem} once by the chosen numerical method. For most discretizations, $\mathcal{N}^2 \lesssim W_{s} \le \mathcal{N}^3$. Finally, $W_{m}$ is the work to evaluate the $X^\calN$-inner product $(f, g)_{X^{\calN}}$ which usually satisfies $\mathcal{N} \lesssim W_{m} \lesssim \mathcal{N}^2$. Using these notations we can present a rough operation count for the three components of the algorithm. \\

\noindent {\bf{Online solve and its preparation:}} The system \eqref{eq:rbstiff} is usually of small size: a set of $N$ linear algebraic equations for $N$ unknowns, with $N \ll \calN$. However, the formation of the stiffness matrix involves $u^{\calN}(\bmu^n)$ for $1 \le n \le N$; direct computation with these quantities requires $\calN$-dependent complexity. It is the affine parameter assumption \eqref{eq:assum_a} that allows us to circumvent complexity in the online stage.  By \eqref{eq:assum_a}, the stiffness matrix for \eqref{eq:rbstiff} can be expressed as
\begin{align}\label{rb_online_system}
  \sum_{m=1}^{N}\sum_{q = 1}^{Q_{a}} \Theta^{q}(\bm{\mu})a^{q}\left( u^{\calN}\left(\bmu^m\right), u^{\calN}\left( \bmu^n\right) \right) u_{Nm}^{\mathcal{N}}(\bm{\mu}) &= f\left(u^{\calN}\left(\bmu^n\right)\right), & n &=1, \ldots, N
\end{align}
During the offline stage, we can precompute the $Q_a$ matrices $a^q\left(u^\calN\left(\bmu^m\right), u^\calN\left(\bmu^n\right)\right)\in \R^{N \times N}$ for $q=1, \ldots, Q_a$ with a cost of order $\calN^2 N^2 Q_a$. 
During the online phase, we need only assemble the reduced stiffness matrix according to \eqref{rb_online_system}, and solve the reduced $N \times N$ system. The total online operation count is thus of order $Q_{a}N^2 + N^3$. \\

\noindent {\bf Greedy sweeping:}  In the offline phase of the algorithm, we repeatedly sweep the training set $\Xi_{\rm train}$ for maximization of the error estimator $\Delta_n(\bmu), \, 1 \le n \le N$. 
The offline cost includes: 
\begin{itemize}[itemsep=0pt]
\item computing the lower bound $\alpha_{LB}^{\mathcal{N}}(\bm{\mu})$. The operation count is  $O(N_{\rm train}W_{\alpha})$, 
\item sweeping the training set by calculating the reduced basis solution and the \textit{a posteriori} error estimate at each location. The operation count $O(N_{\rm train}Q^2_{a}N^{3}_{RB})$
\item solving system \eqref{eq:update_problem} $N$ times. The total operation count is $O(N W_{s})$.
\end{itemize}

\noindent {\bf Error estimator calculations:} 
With a cost of order $Q_a N \calN$ in the offline stage, we can calculate functions $C$ and $\mathcal{L}_{m}^{q}$, $1\le m \le N, 1 \le q \le Q_{a}$ both defined by
\begin{equation}\label{error_es_problem}
\begin{cases}
  (\mathcal{C}, v) = f(v)_{X^{\calN}} \quad \forall v \in X^{\mathcal{N}} \\
  (\mathcal{L}_{m}^{q}, v)_{X^{\calN}} = -a^{q}(u^{\calN}\left(\bmu^m\right), v) \quad \forall v \in X^{\mathcal{N}}.
\end{cases}
\end{equation}
Here, we assume that the $X$-inner product can be ``inverted'' with cost of order $\calN$, i.e. that the mass matrix is block diagonal. The availability of $\calC$ and $\calL_m^q$ facilitates an Offline-Online decomposition of the term $\lVert r_N(\cdot; \bm{\mu})\rVert_{(X^{\mathcal{N}})'}$ in the error estimate \eqref{eq:error_estimator_en} due to that its square can be written as
\begin{equation}
\label{eq:error_es_quan}
(\mathcal{C}, \mathcal{C})_{X^{\calN}} + 2\sum_{q = 1}^{Q_{a}}\sum_{m = 1}^{N}\Theta^{q}(\bm{\mu})u^{\mathcal{N}}_{Nm}(\bm{\mu})(\mathcal{C}, \mathcal{L}_{m}^{q})_{X} + \sum_{q = 1}^{Q_{a}}\sum_{m = 1}^{N}\Theta^{q}(\bm{\mu})u^{\mathcal{N}}_{Nm} \left\{ \sum_{q' = 1}^{Q_{a}}\sum_{m' = 1}^{N}\Theta^{q'}(\bm{\mu})u^{\mathcal{N}}_{Nm'}(\mathcal{L}_{m}^{q}, \mathcal{L}_{m'}^{q'})_{X^{\calN}}\right\}.
\end{equation}

Therefore, in the offline stage we should calculate and store $(\mathcal{C}, \mathcal{C})_{X^{\calN}}, (\mathcal{C},\mathcal{L}_{m}^{q})_{X^{\calN}}, (\mathcal{L}_{m}^{q}, \mathcal{L}_{m'}^{q'})_{X^{\calN}}, \,\, 1 \le m, m' \le N_{\rm RB}, 1 \le q, q' \le Q_{a}$. 
This cost is of the order $Q_a^2 N^3 W_m$.
During the online stage, given any parameter $\bm{\mu}$, we only need to evaluate $\Theta^{q}(\bm{\mu}), 1 \le q \le Q, u^{\mathcal{N}}_{Nm}(\bm{\mu}), 1 \le m \le N$, and compute the sum \eqref{eq:error_es_quan}. Thus, the online operation count for each $\bm{\mu}$ is $O(Q^2_{a}N^3)$.

\noindent {\bf Summary}, the total offline portion of the algorithm has complexity of the order
\begin{align*}
  \overbrace{\calN^2 N^2 Q_a}^{\mbox{Reduced solve preparation}} + \overbrace{N_{\rm train} W_\alpha + N_{\rm train} (Q_a^2 N^3 + N^4)+ N W_s}^{\mbox{Greedy sweeping}} + \overbrace{Q_a^2 N^3 W_m}^{\mbox{Estimator preparation}}.
\end{align*}
The total online cost including the error certification is of order $Q^2_{a}N^3$.

\subsection{Offline-Online decomposition for the approximate CDM-RBM Gramian $\widetilde{G}$}
\label{sec:app-offline-online-CDG}

The entries of the matrix $\widetilde{G}$ defined in \eqref{eq:error_matrix} can be efficiently computed assuming that we can compute the approximate errors $\left\{\widetilde{e}(\bmu): \bmu \in \Xi_{\rm train}\right\}$ in an offline-online fashion. To accomplish this, note that 
\begin{equation*}
\begin{aligned}
  \widetilde{e}(\bmu)  & =  \bigg(\sum_{m = 1}^{Q}u^{\mathcal{N}}_{Nm}(\bm{\mu}) \mathbb{A}^{-1}_{\mathcal{N}}(\bm{\mu}^{m})\bigg) \bigg(f^{\mathcal{N}} - \mathbb{A}_{\mathcal{N}}(\bm{\mu})\left(\sum^{N}_{m=1}u^{\mathcal{N}}_{Nm}(\bm{\mu}) u^{\calN}\left(\bmu^m\right) \right)\bigg)\\
                       & =\sum_{m = 1}^{Q} u^{\mathcal{N}}_{Nm}(\bm{\mu}) \bigg( \mathbb{A}^{-1}_{\mathcal{N}}(\bm{\mu}^{m})f^{\mathcal{N}}\bigg) - \sum_{m = 1}^{Q}\sum_{m' = 1}^{N} u^{\mathcal{N}}_{Nm}(\bm{\mu})u^{\mathcal{N}}_{Nm'}(\bm{\mu}) \bigg(\mathbb{A}^{-1}_{\mathcal{N}}(\bm{\mu}^{m})\mathbb{A}_{\mathcal{N}}(\bm{\mu})\left(u^{\calN}\left(\bmu^{m'}\right)\right)\bigg)\\
                       &  = \sum_{m = 1}^{Q} u^{\mathcal{N}}_{Nm}(\bm{\mu}) \bigg( \mathbb{A}^{-1}_{\mathcal{N}}(\bm{\mu}^{m})f^{\mathcal{N}}\bigg) - \sum_{m = 1}^{Q}\sum_{m' = 1}^{N} \sum_{k = 1}^{Q_{a}} \theta^a_{k}(\bm{\mu})u^{\mathcal{N}}_{Nm}(\bm{\mu})u^{\mathcal{N}}_{Nm'}(\bm{\mu}) \bigg(\mathbb{A}^{-1}_{\mathcal{N}}(\bm{\mu}^{m})A_{k}( u^{\calN}\left(\bmu^{m'}\right))\bigg),
\end{aligned}
\end{equation*}
Therefore, we can split this computation into offline and online components as follows:
\begin{itemize}
  \item {\bf Offline:} Calculate  $\mathbb{A}^{-1}_{\mathcal{N}}(\bm{\mu}^{m})f^{\mathcal{N}}$ and $\mathbb{A}^{-1}_{\mathcal{N}}(\bm{\mu}^{m})A_{k}\left(u^{\calN}\left(\bmu^{m'}\right)\right)$ for $1 \le m' \le N,~1 \le m \le Q ~,~ 1 \le k \le Q_{a},~ {Q \le N}$, with complexity $O(\mathcal{N}^2Q N Q_{a})$.
\item {\bf Online:} Evaluate the coefficients $u^{\mathcal{N}}_{Nm}(\bm{\mu})$ and $\theta^a_{k}(\bm{\mu})u^{\mathcal{N}}_{Nm}(\bm{\mu})u^{\mathcal{N}}_{Nm'}(\bm{\mu})$ and form $\widetilde{e}(\bmu)$. The online computation has complexity $O(\mathcal{N}Q N Q_{a})$.
\end{itemize}

\begin{figure}
\begin{center}
  \resizebox{\textwidth}{!}{
    \includegraphics[width=\textwidth]{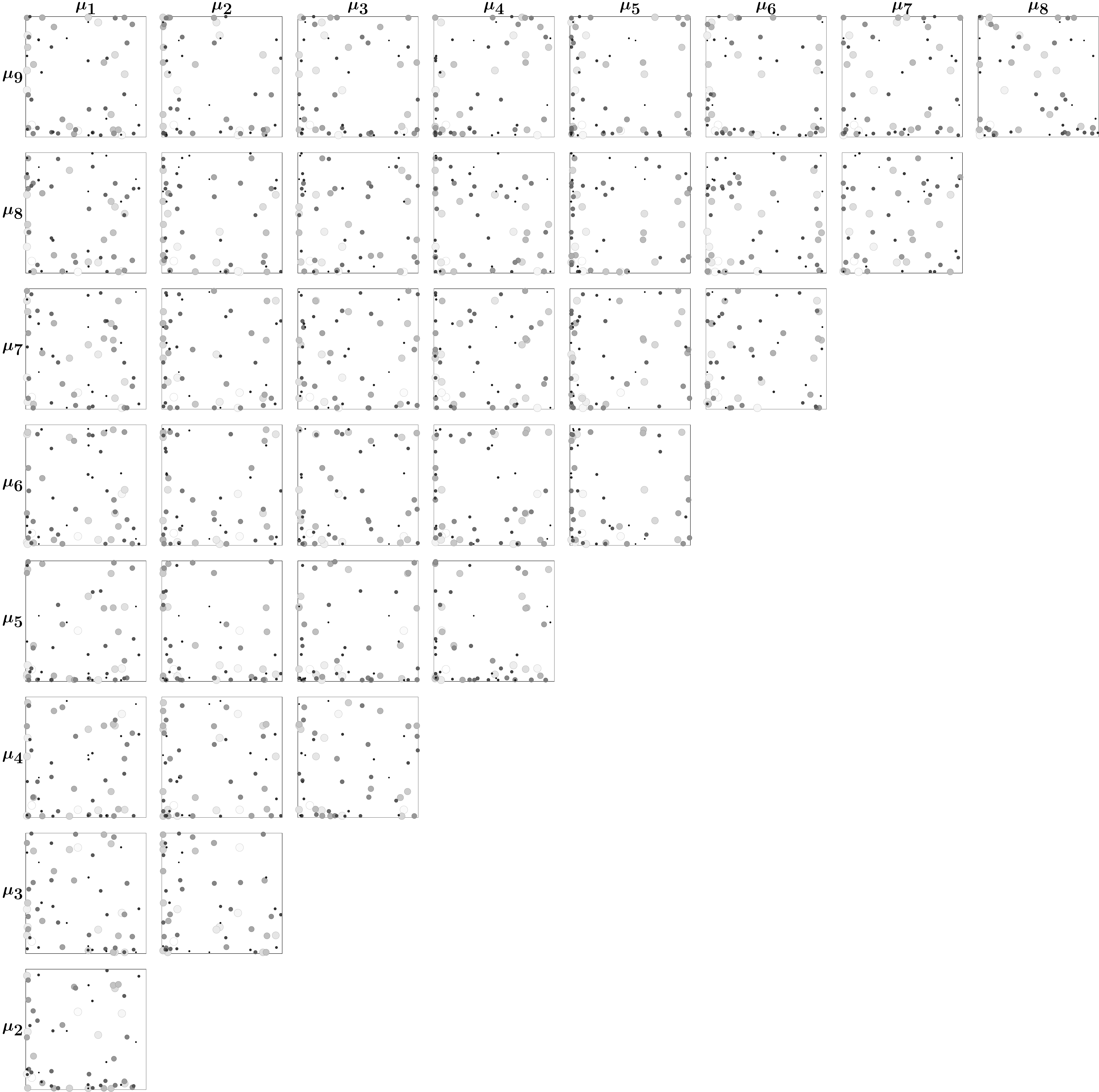}
  }
\end{center}
\caption{Parameter selections for the test problem 2 with a 9-dimensional parameter. Two-dimensional scatter plots of the parameter selection are shown corresponding to the classical RBM algorithm. The size and transparency of the markers indicate the ordering of the sequence: points earlier in the sequence have a larger radius and are more transparent, points later in the sequence are smaller and more opaque. Both the horizontal and vertical axes in each plot range over the interval $[0.1, 10]$.}
\label{fig:bivariate-marginals-problem2-rb}
\end{figure}

\begin{figure}
\begin{center}
  \resizebox{\textwidth}{!}{
    \includegraphics[width=\textwidth]{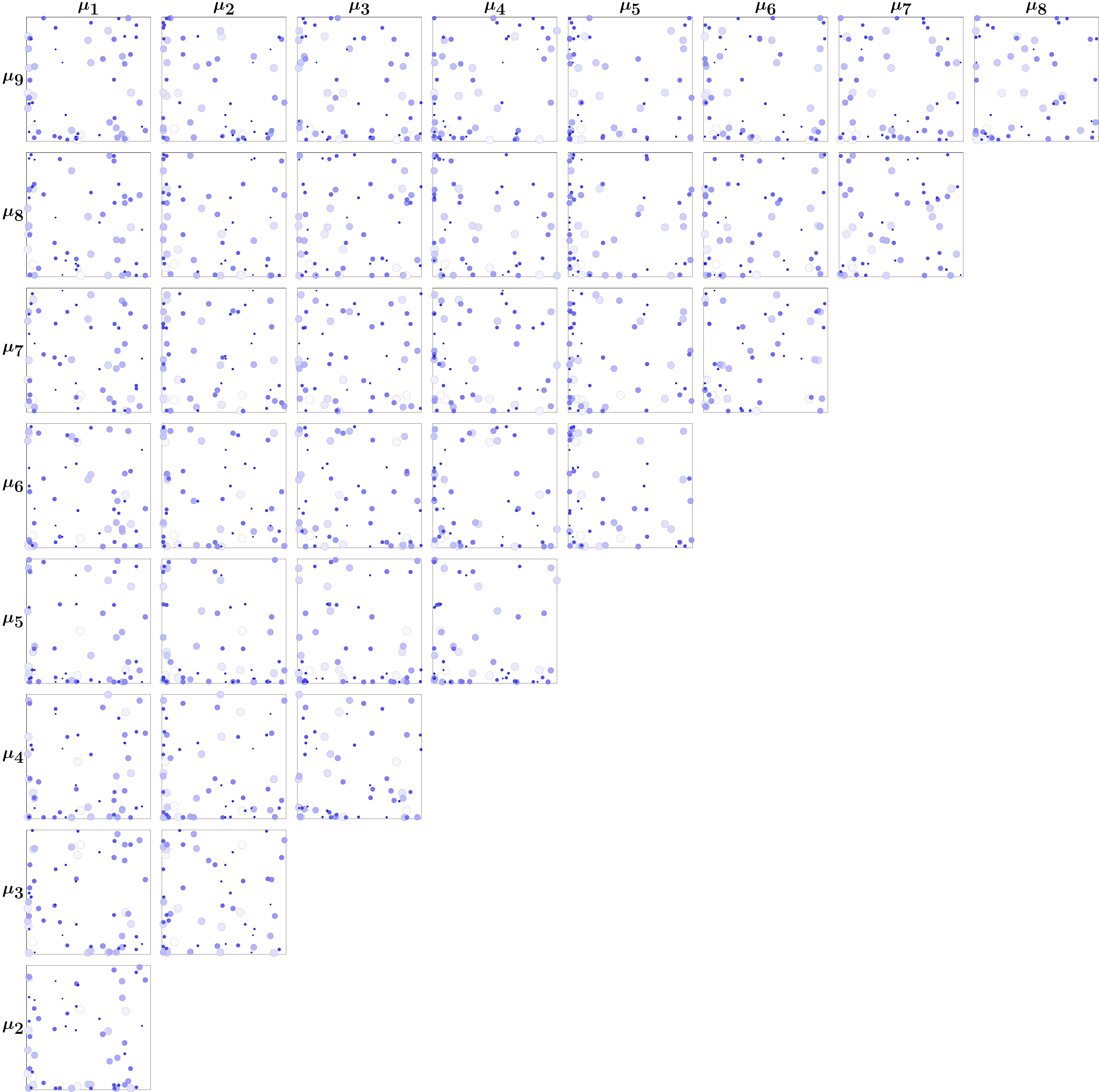}
  }
\end{center}
\caption{Parameter selections for the test problem 2 with a 9-dimensional parameter. Two-dimensional scatter plots of the parameter selection are shown corresponding to the SMM-RBM algorithm. The size and transparency of the markers indicate the ordering of the sequence: points earlier in the sequence have a larger radius and are more transparent, points later in the sequence are smaller and more opaque. Both the horizontal and vertical axes in each plot range over the interval $[0.1, 10]$.}
\label{fig:bivariate-marginals-problem2-smm}
\end{figure}

\begin{figure}
\begin{center}
  \resizebox{\textwidth}{!}{
    \includegraphics[width=\textwidth]{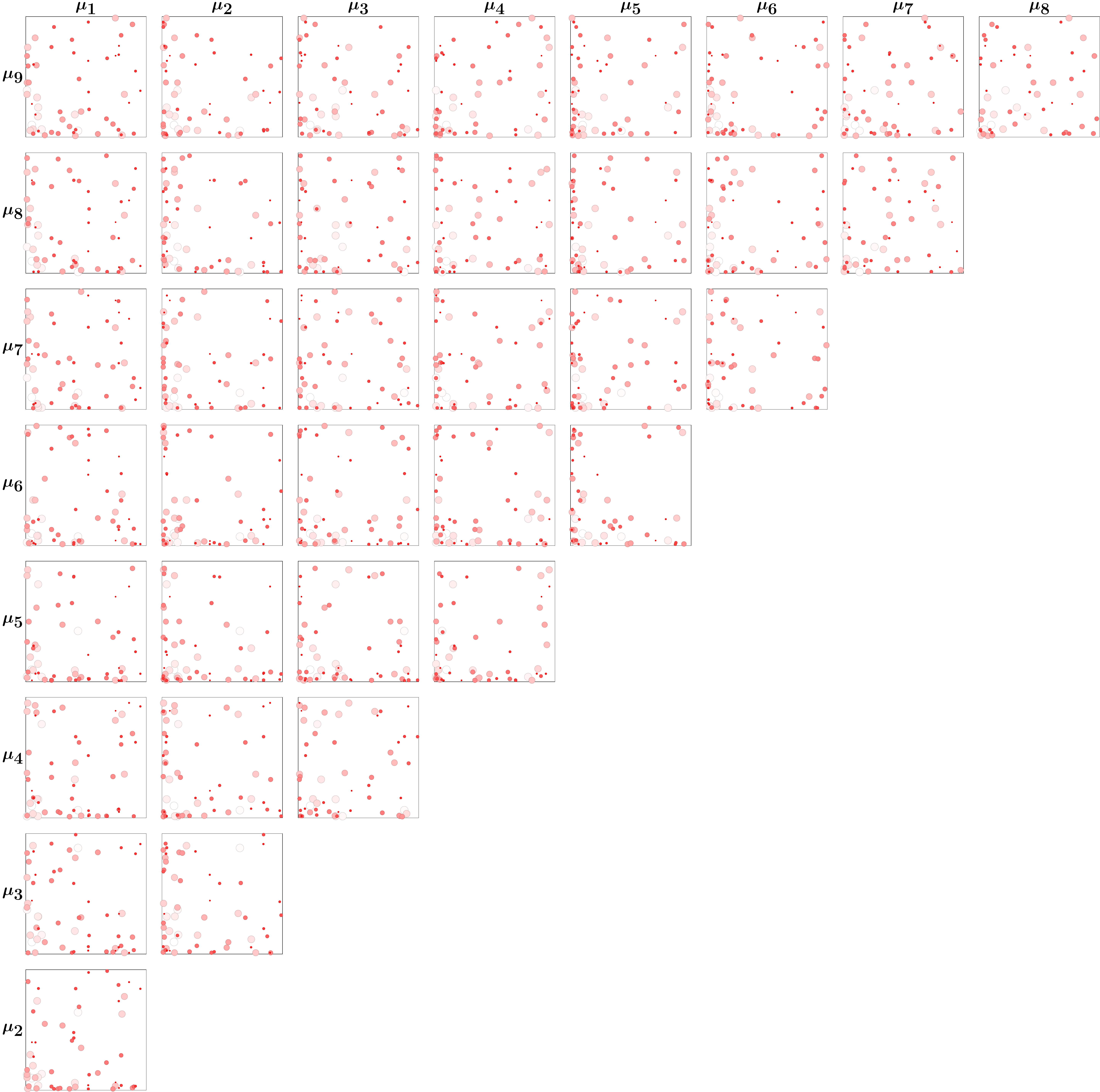}
  }
\end{center}
\caption{Parameter selections for the test problem 2 with a 9-dimensional parameter. Two-dimensional scatter plots of the parameter selection are shown corresponding to the CDM-RBM algorithm. The size and transparency of the markers indicate the ordering of the sequence: points earlier in the sequence have a larger radius and are more transparent, points later in the sequence are smaller and more opaque. Both the horizontal and vertical axes in each plot range over the interval $[0.1, 10]$.}
\label{fig:bivariate-marginals-problem2-cdm}
\end{figure}

\bibliographystyle{abbrv}

\bibliography{multibib}

\end{document}